\documentclass[a4paper, 12pt]{article}
\usepackage{enumerate, theorem}
\usepackage{amsmath, amsfonts, amssymb}
\usepackage[height=22.5cm, width=15cm]{geometry}
\usepackage[all]{xy}
\usepackage{mathrsfs, yfonts}

\DeclareMathOperator{\C}{\mathbb{C}}

\newcommand{\A}{\tilde{\mathcal{A}}}

\newcommand{\parag}[1]{\paragraph{\sc{#1.}}}

\newtheorem{thm}{Theorem}[subsection]
\newtheorem{defn}[thm]{Definition}
\newtheorem{cor}[thm]{Corollary}
\newtheorem{prop}[thm]{Proposition}
\newtheorem{lemma}[thm]{Lemma}

\begin{document}
\title{Asymptotics of a vanishing period : \\
characterization of semi-simplicity}

\author{Daniel Barlet\footnote{Barlet Daniel, Institut Elie Cartan UMR 7502  \newline
Universit\'e de Lorraine, Institut Elie Cartan de Lorraine, UMR 7502, Vandoeuvre-les-Nancy, F-54506, France \newline
CNRS,  Institut Elie Cartan de Lorraine, UMR 7502, Vandoeuvre-les-Nancy, F-54506, France \newline
  Institut Universitaire de France, \newline
e-mail : Daniel.Barlet@univ-lorraine.fr}.}

\date{27/1/13.}

\maketitle

\section*{Abstract}
In this paper we introduce the word {\em fresco} to denote a monogenic geometric (a,b)-module. This "basic object"  (generalized Brieskorn module with one generator)  corresponds to the formal germ of the minimal filtered (regular) differential equation. Such an equation is satisfied by a relative de Rham cohomology class at a critical value of a holomorphic function on a smooth complex manifold. In [B.09]  the first structure theorems are proved. Then in [B.10] we introduced the notion of {\em theme} which corresponds in the \ $[\lambda]-$primitive case  to frescos having a unique Jordan-H{\"o}lder sequence (a unique Jordan block for the monodromy). Themes correspond  to asymptotic expansion of a given vanishing period, so to an image of a fresco in the module of  asymptotic expansions. For a fixed relative de Rham cohomology class (for instance given by  a smooth differential form $d-$closed and $df-$closed) each choice of a vanishing cycle in the spectral eigenspace of the monodromy for the eigenvalue \ $exp(2i\pi.\lambda)$ \ produces a \ $[\lambda]-$primitive theme, which is a quotient of the fresco associated to the given relative de Rham class  itself. \\
We show that for any  fresco there exists an {\em unique} Jordan-H{\"o}lder sequence, called the {\em principal J-H. sequence}, with corresponding quotients giving the opposite of  the roots of the Bernstein polynomial in  increasing order. We  study  the semi-simple part of a given fresco and we  characterize the semi-simplicity of a fresco  by the fact  for any given order on the roots of its  Bernstein polynomial we may find a J-H. sequence making them appear with this order. Then we construct a numerical invariant, called the \ $\beta-$invariant, and we show that it produces numerical criteria in order to give a necessary and sufficient condition on a fresco to be semi-simple. We show that these numerical invariants define a natural algebraic stratification on the set of isomorphism classes of fresco with given fundamental invariants (or equivalently with given roots of the Bernstein polynomial).

\parag{AMS Classification} 32 S 25, 32 S 40, 32 S 50.

\parag{Key words} Fresco,  theme,  (a,b)-module, Brieskorn module, asymptotic expansion, vanishing 
period, Gauss-Manin connection, filtered differential equation, Bernstein polynomial.

\tableofcontents

\section*{Introduction}
Let \ $f : X \to D$ \ be an holomorphic function on a connected  complex manifold. Assume that \ $\{df = 0 \} \subset \{ f = 0 \}: = X_0 $. We consider \ $X$ \ as a degenerating family of complex manifolds parametrized by \ $D^* : = D \setminus \{0\}$ \ with a singular member  \ $X_0$ \ at the origin of \ $D$. Let \ $\omega$ \ be a smooth \ $(p+1)-$differential form on \ $X$ \ satisfying \ $d\omega = 0 = df\wedge \omega$. Then in many interesting cases (see for instance [Br.70], [M.74], [S.89] for the isolated singularity case, [B.I], [B.II], [B.III] for the case of a function with 1-dimensional singular set and  [B.12]  for the general proper case) the relative family of de Rham cohomology classes induced on the fibers \ $(X_s)_{ s \in D^*}$ \ of \ $f$ \ by  \ $\omega\big/df$ \ is solution of a minimal filtered differential equation defined from the Gauss-Manin connection of \ $f$. This object, called a {\bf fresco}, is a monogenic\footnote{i.e. a left \ $\A-$module generated by one element.} regular (a,b)-module satisfying an extra condition, called "geometric", which encodes simultaneously the regularity at \ $0$ \ of the Gauss-Manin connection, the monodromy theorem and B. Malgrange's positivity theorem (see [M.75], the result [K.76] on the roots of the Bernstein polynomial or [Bj.93]).

\smallskip

After some easy preliminaries, we prove in section 1 that a fresco \ $E$ \  admits an {\em unique} Jordan-H{\"o}lder sequence, called the {\em principal} J-H. sequence, in which the opposite of the roots of the Bernstein polynomial of \ $E$ \ appears in the increasing order\footnote{We fix an order on \ $\mathbb{Q}$ \ deduced from a total order on \ $\mathbb{Q}\big/\mathbb{Z}$ \  and compatible with the usual order on each class modulo \ $\mathbb{Z}$.}. This uniqueness result is important because it implies, for instance, that the isomorphism class of each quotient of two terms of the principal J-H. sequence only depends on the isomorphism class of \ $E$.

\smallskip

In section 2 we  study semi-simple regular (a,b)-modules and the corresponding semi-simple filtration. In the case of a   fresco we prove that semi-simplicity is characterized by the fact that we may find a J-H. sequence in which the opposite of the roots of the Bernstein polynomial appears in strictly decreasing order (but also in any given order). We conclude this section by an example of fresco such that the semi-simple part of \ $E$ \ and of its dual  \ $E^{*}$ \ have not the same rank.

\smallskip

In section 3 we answer to the following question : How to recognize from the principal J-H. sequence of a \ $[\lambda]-$primitive fresco if it is semi-simple? Of course the answer is easy if some rank \ $2$ \ sub-quotient theme appears from this sequence. But when it is not the case, such a rank \ $2$ \ sub-quotient theme may appear after ''commuting'' some terms in the J-H. sequence. The simplest example is when
$$ E : = \A\big/\A.(a - \lambda_1.b)(1 + \alpha.b^{p_1+p_2})^{-1}.(a - \lambda_2.b).(a - \lambda_3.b) $$
with \ $\lambda_{i+1} = \lambda_i + p_i -1$ \ for \ $i = 1,2$ \ with \ $p_i \in \mathbb{N}^*,  i = 1,2$ \ and \ $\alpha \in \C^*$. Using the identity in \ $\A$
 $$(a - \lambda_2.b).(a - \lambda_3.b) = (a - (\lambda_3+1).b).(a - (\lambda_2-1).b) $$
 it is easy to see that \ $(a - (\lambda_2-1).b).[1]$ \ generates a (normal) rank \ $2$ \ theme in \ $E$, because we assume \ $\alpha \not= 0$.
 
 \smallskip
 
We solve this question introducing for a \ $[\lambda]-$primitive rank \ $k \geq 3$ \  fresco \ $E$ \ such that \ $F_{k-1}(E)$ \ and \ $E\big/F_1(E)$ \ are semi-simple, where \ $(F_j(E))_{j \in [1,k]}$ \ is the principal J-H. sequence of \ $E$, the \ $\beta-$invariant \ $\beta(E)$. This complex number only depends on the isomorphism class of \ $E$ \ and is zero if and only \ $E$ \ is semi-simple. We also prove that when \ $\beta(E)$ \ is not zero it determines the isomorphism class of any normal rank \ $2$ \ sub-theme of \ $E$. This easily implies that the polynomial  \ $\beta$ \ is quasi-invariant by change of variables.
We use this \ $\beta-$invariant as a  main tool to describe, in the set \ $\mathcal{F}(\lambda_1, \dots, \lambda_k)$ \ of isomorphism classes of rank \ $k$ \  frescos with fundamental invariants \ $\lambda_1, \dots, \lambda_k$, an {\em algebraic}\footnote{notion defined in the section 3.} stratification ending in the subset of semi-simple frescos. This stratification is also invariant by any change of variable.

\section{Preliminaries.}

\subsection{Definitions and characterization as \ $\A-$modules.}

We are interested in "standard " formal asymptotic expansions of the following type
$$ \sum_{q = 1}^k \sum_{j=0}^N \ \C[[s]].s^{\lambda_q-1}.(Log\, s)^j $$
where \ $\lambda_1, \dots, \lambda_k$ \ are positive rational numbers, and in fact in vector valued such expansions. The two basic operations on such expansions are
\begin{itemize}
\item the multiplication by \ $s$ \ that  we shall denote \ $a$,
\item and the primitive in \ $s$ \ without constant that we shall denote \ $b$.
\end{itemize}
This leads to consider on the set of such expansions a left module structure on the  non commutative \ $\C-$algebra 
$$ \A : = \{\  \sum_{\nu = 0}^{\infty} \ P_{\nu}(a).b^{\nu},\quad {\rm where} \quad P_{\nu} \in \C[x] \  \}  $$
defined by  the following conditions
\begin{itemize}
\item The commutation relation \ $a.b - b.a = b^2$ \ which is the translation of the Leibnitz rule ;
\item The continuity for the \ $b-$adic filtration of \ $\A$ \ of the left and right multiplications by \ $a$.
\end{itemize}
Define now for \ $\lambda$ \ in \ $\mathbb{Q}\, \cap\, ]0,1]$ \ and \ $N \in \mathbb{N}$ \  the left \ $\A-$module
$$ \Xi^{(N)}_{\lambda} : = \oplus_{j=0}^N \ \C[[s]].s^{\lambda-1}.(Log\, s)^j  = \oplus_{j=0}^N \ \C[[a]].s^{\lambda-1}.(Log\, s)^j = \oplus_{j=0}^N \ \C[[b]].s^{\lambda-1}.(Log\, s)^j .$$
Of course we let \ $a$ \ and \ $b$ \ act on \ $\Xi_{\lambda}^{(N)}$ \ as explained above.\\
Define also, when \ $\Lambda$ \ is a finite  subset in \ $\mathbb{Q}\, \cap\, ]0,1]$, the \ $\A-$module
$$ \Xi_{\Lambda}^{(N)} : = \oplus_{\lambda \in \Lambda}\  \Xi_{\lambda}^{(N)} .$$
More generally, if \ $V$ \ is a finite dimensional complex vector space we shall put a structure of left \ $\A-$module on \ $\Xi_{\Lambda}^{(N)}\otimes_{\C} V$ \ with the following rules :
$$ a.(\varphi \otimes v) = (a.\varphi)\otimes v \quad {\rm and} \quad b.(\varphi \otimes v) = (b.\varphi)\otimes v$$
for any \ $\varphi \in \Xi_{\Lambda}^{(N)}$ \ and any \ $v \in V$. It will be convenient to denote  \ $\Xi_{\lambda}$ \ the \ $\A-$module \ $\sum_{N \in \mathbb{N}}\ \Xi^{(N)}_{\lambda}$. \\

\begin{defn}\label{fresco}
A \ $\A-$module is call a {\bf fresco} when it is isomorphic to a submodule \ $\A.\varphi \subset \Xi_{\Lambda}^{(N)}\otimes V$ \ where \ $\varphi$ \ is any element in \ $\Xi_{\Lambda}^{(N)}\otimes V$, for some choice of \ $\Lambda, N$ \ and \ $V$ \ as above.\\
A fresco is a {\bf theme} when we may choose \ $V : = \C$ in the preceeding choice.
\end{defn}
For the motivation of these definition see the fundamental example below.

\bigskip

Now the characterization of frescos among all left \ $\A-$modules is not so obvious. The following theorem is proved in [B.09].

\begin{thm}\label{charact. frescos}
A left \ $\A-$module \ $E$ \  is a fresco if and only if it is a geometric (a,b)-module which is generated (as a \ $\A-$module) by one element. Moreover the annihilator in \ $\A$ \ of  a generator of \ $E$ \ is a left ideal of the form \ $\A.P$ \ where \ $P$ \ may be written as follows :
$$ P = (a - \lambda_1.b).S_1^{-1}.(a - \lambda_2.b).S_2^{-1} \dots (a - \lambda_k.b).S_k^{-1}, \quad k : =  \dim_{\C}(E\big/b.E) $$
where \ $\lambda_j$ \ are rational numbers such that \ $\lambda_j + j > k$ \ for \ $j \in [1,k]$ \ and where \ $S_1, \dots , S_k$ \ are invertible elements in the sub-algebra \ $\C[[b]]$ \ of \ $\A$.\\
Conversely, for such a \ $P \in \A$ \ the left \ $\A-$module \ $E : = \A\big/\A.P$ \ is a fresco and it is a free rank \ $k$ \ module on \ $\C[[b]]$. 
\end{thm}

Let me recall briefly for the convenience of the reader the definitions of the notions involved in the previous statement.
\begin{itemize}
\item A {\bf (a,b)-module} \ $E$ \ is a free finite rank \ $\C[[b]]$ \ module endowed with an \ $\C-$linear endomorphism \ $a$ \ such that \ $a.S = S.a + b^2.S'$ \ for \ $S \in \C[[b]]$;  or, in an equivalent way, a \ $\A-$module which is free and finite type over the subalgebra \ $\C[[b]] \subset \A$. \\
 It has a {\bf simple pole} when \ $a$ \ satisfies \ $a.E \subset b.E$. In this case the {\bf Bernstein polynomial} \ $B_E$ \  of \ $E$ \ is defined as the {\em minimal} polynomial of \ $-b^{-1}.a$ \ acting on \ $E\big/b.E$.
\item A (a,b)-module \ $E$ \  is {\bf regular}  when it may be embedded in a simple pole (a,b)-module. In this case there is a minimal such embedding which is the inclusion of \ $E$ \ in its saturation \ $E^{\sharp}$ \ by \ $b^{-1}.a$. The {\bf Bernstein polynomial} of a regular (a,b)-module is, by definition, the Bernstein polynomial of its saturation \ $E^{\sharp}$.
\item A regular (a,b)-module \ $E$ \ is called {\bf geometric} when all roots of its Bernstein polynomial are rational and strictly negative (compare with [K.76]).
\end{itemize}

Note that the formal completion in \ $b$ \ of the Brieskorn module of a function with an isolated singularity is a geometric (a,b)-module (see [Br.70], [M.74]). In  [B.I], [B.II] and [B.III] we have used this notion to study  holomorphic germ in \ $\C^{n+1}$ \ with a \ $1-$dimensional singularity. The result in [B.12] shows that this structure appears in a rather systematic way in the study of the Gauss-Manin connection of a proper  holomorphic function on a complex manifold.

\parag{Fundamental example} Let \ $E$ \ be the formal completion of the Brieskorn module of a holomorphic germ \ $f :( \C^{n+1},0) \to (\C,0)$ \ with an isolated singularity. Then for any  \ $(n+1)-$holomorphic form \ $\omega$ \ we have a fresco \ $\A.[\omega] \subset E$. Then for each \ $\gamma$,  a $n-$homology class of the Milnor fiber of \ $f$, we have a \ $\A-$linear map \ $ \varphi_{\gamma} : E \to \Xi$ \ given by \ $\varphi_{\gamma}(\omega) : = \int_{\gamma_{s}} \omega\big/df $, where \ $\Xi$ \ is the \ $\A-$module of asymptotic expansions and \ $(\gamma_{s})_{s \in D^{*}}$ \ the horizontal multivalued family of $n-$cycles in the fibers of \ $f$ \ associated to \ $\gamma$. Then for a given class \ $[\omega] \in E$ \ the image by \ $\varphi_{\gamma}$ \ of  the fresco \ $\A.[\omega]$ \ is a theme. If we choose \ $\gamma$ \ in the spectral subspace of the the monodromy for the eigenvalue \ $exp(2i\pi\lambda)$, this theme will be \ $[\lambda]-$primitive (see the definition above). In this case the rank of such a theme is bounded by the nilpotency order of  the monodromy acting on \ $\gamma$.\hfill $\square$\\

\begin{defn}\label{normal}A submodule \ $F $ \ in \ $E$ \ is {\bf normal} when \ $F \cap b.E = b.F$.
\end{defn}

For any sub-module \ $F$ \ of a (a,b)-module \ $E$ \ there exists a minimal normal sub-module \ $\tilde{F}$ \ of \ $E$ \ containing \ $F$. We shall call it the {\bf normalization} of \ $F$. It is easy to see that \ $\tilde{F}$ \ is the pull-back by the quotient map \ $E \to E\big/F$ \ of the \ $b-$torsion of \ $E\big/F$. Note that \ $F$ \ is always a finite codimensional complex vector space of its normalization. In particular \ $F$ \ and \ $\tilde{F}$ \  have the same rank as \ $\C[[b]]-$modules.

\smallskip

 When \ $F$ \ is normal the quotient \ $E\big/F$ \ is again a (a,b)-module. Note that a sub-module of a regular (resp. geometric) (a,b)-module is regular (resp. geometric) and when \ $F$ \ is normal \ $E\big/F$ \ is also regular (resp. geometric).

\begin{lemma}
Let \ $E$ \ be a fresco and \ $F$ \ be a normal sub-module in \ $E$. Then \ $F$ \ is a fresco and also the quotient \ $E\big/F$.
\end{lemma}

\parag{proof} The only point to prove, as we already know that \ $F$ \ and \ $E\big/F$ \ are geometric (a,b)-modules thanks to [B.09], is the fact that \ $F$ \ and \ $E\big/F$ \ are generated as \ $\A-$modules by one element. This is obvious for \ $E\big/F$, but not for \ $F$. We shall use the theorem \ref{charact. frescos} for \ $E\big/F$ \ to prove that \ $F$ \ is generated by one element. Let \ $e$ \ be a generator of \ $E$ \ and let \ $P$ \ as in the theorem \ref{charact. frescos} which generates the annihilator ideal of the image of \ $e$ \ in \ $E\big/F$. Then \ $P.e$ \ is in \ $F$. We shall prove that \ $P.e$ generates \ $F$ \ as a \ $\A-$module. Let \ $y$ \ be an element in \ $F$ \ and write \ $y = u.e$ \ where \ $u$ \ is in \ $\A$. As \ $P$ \ is, up to an invertible element in \ $\C[[b]]$, a monic polynomial in \ $a$ \ with coefficients in \ $\C[[b]]$, we may write \ $u = Q.P + R$ \ where \ $Q$ \ and \ $R$ \ are in \ $\A$ \ and \ $R$ \ is a polynomial in \ $a$ \ with coefficient in \ $\C[[b]]$ \ of degree \ $r < deg(P) = rank(E\big/F)$. Now, as \ $y$ \ is in \ $F$,  the image in \ $E\big/F$ \ of \ $u.e$ \ is \ $0$, and this implies that \ $R$ \ annihilates the image of \ $e$ \ in \ $E\big/F$. So \ $R$ \ lies in \ $\A.P$ \ and so \ $R = 0$. Then we have \ $u = Q.P$ \ and \ $y = Q.P.e $ \ proving our claim. $\hfill \blacksquare$\\

The following lemma is useful to recognize a generator of a fresco.

\begin{lemma}\label{gen.}
Let \ $E$ \ be a non zero geometric (a,b)-module. Then \ $E$ \ is a fresco if and only if  the vector space \ $E\big/a.E+b.E$ \ has dimension \ $1$. In this case any element in \ $E \setminus a.E + b.E$ \ is a generator of \ $E$ \ as a \ $\A-$module.
\end{lemma}

\parag{proof}  If \ $E$ \ is a rank \ $k \geq 1$ \ fresco then \ $a$ \ acts on the \ $k-$dimensional vector space \ $E\big/.b.E$ \ as a principal nilpotent endomorphism (i. e. its minimal polynomial is \ $ P(x) = x^{k}$). So for any \ $e \not\in a.E + b.E$ \ the set \ $(e, a.e, \dots, a^{k-1}.e)$, induces a basis of \ $E\big/b.E$ \ and then, it is a \ $\C[[b]]-$basis of \ $E$. Such an element \ $e$ \ is a generator of \ $E$ \ as a left \ $\A-$module.\\
Conversely, assume that, for the rank \ $k$ \ geometric (a,b)-module \ $E$, the vector space \ $E\big/a.E + b.E$ \ is \ $1-$dimensional;  then \ $a$ \ acts as a principal nilpotent endomorphism on \ $E\big/b.E$, and  for any \ $e \in E \setminus a.E + b.E$ \  the set \ $(e, a.e, \dots, a^{k-1}.e)$ \ is a \ $\C[[b]]-$basis of \ $E$. So such an element generates \ $E$ \ as a left \ $\A-$module.$\hfill\blacksquare$\\

In the case of a fresco the Bernstein polynomial is more easy to describe, thanks to the following proposition proved in [B.09].

\begin{prop}\label{Bernst. fresco}
Soit \ $E = \A\big/\A.P$ \ be a rank \ $k$ \ fresco as described in the previous theorem. The Bernstein polynomial of \ $E$ \ is the {\bf characteristic} polynomial of \ $-b^{-1}.a$ \ acting on \ $E^{\sharp}\big/b.E^{\sharp}$. And  the {\bf Bernstein element}  \ $P_E$ \ of \ $E$, which is the element in \ $\A$ \ defined by the Bernstein polynomial \ $B_E$ \ of \ $E$ \ by the following formula
$$ P_E : = (-b)^k.B_E(-b^{-1}.a) $$
is equal to \ $(a - \lambda_1.b) \dots (a - \lambda_k.b)$ \ for any such a choice of presentation of \ $E$ \ as in the theorem \ref{charact. frescos}
\end{prop}

As an easy consequence, the \ $k$ \  roots of the Bernstein polynomial of \ $E$ \ are the opposite of the numbers \ $\lambda_1+1-k, \dots, \lambda_k+k-k$. So the Bernstein polynomial is readable on the element \ $P$ : the initial form of \ $P$ \ in (a,b) has degree \ $k$ \ and  is the Bernstein element\ $P_E$ \ of \ $E$. 

\parag{Remark} If we have an exact sequence of (a,b)-modules
$$ 0 \to F \to E \to G \to 0 $$
where \ $E$ \ is a fresco, then \ $F$ \ and \ $G$ \ are frescos and the Bernstein elements satisfy the equality \ $P_E = P_F.P_G$ \ in the algebra \ $\A$ \ (see [B.09]  proposition 3.4.4.). $\hfill \square$

\parag{Example} The \ $\A-$module \ $\Xi_{\lambda}^{(N)}$ \ is a simple pole (a,b)-module with rank \ $N+1$. Its Bernstein polynomial is equal to \ $(x+\lambda)^{N+1}$.\\
The theme \ $\A.\varphi \subset \Xi^{(N)}_{\lambda}$ \ where \ $\varphi = s^{\lambda-1}.(Log\,s)^N$ \ has rank \ $N+1$ \ and its Bernstein element is \ $(a - (\lambda+N).b)(a - (\lambda+N-1).b) \dots (a - \lambda.b)$. Its saturation by \ $b^{-1}.a$ \  is \ $\Xi^{(N)}_{\lambda}$.$\hfill \square$\\

The dual of a regular (a,b)-module is regular, but duality does not preserve the property of being geometric because duality changes the sign of the roots of the Bernstein polynomial. As duality preserves regularity (see [B.95]), to find again a geometric (a,b)-module it is sufficient to make the tensor product  by a rank \ $1$ \ (a,b)-module \ $E_{\delta}$ \ for \ $\delta$ \ a large enough rational number\footnote{to tensor by \ $E_{\delta}$ \ is the same that to replace \ $a$ \ by \ $a + \delta.b$ ; see [B.I] for the general  definition of the tensor product.}. The next lemma states that this "twist duality" preserves the notion of fresco.

\begin{lemma}\label{twisted duality}
Let \ $E$ \ be a fresco and let \ $\delta$ \ be a rational number such that \ $E^*\otimes E_{\delta}$ \ is geometric. Then \ $E^*\otimes E_{\delta}$ \ is a fresco. 
\end{lemma}

The proof is obvious. $\hfill \blacksquare$\\

The following definition is useful when we want to consider only the part of asymptotic expansions corresponding to prescribe eigenvalues of the monodromy.

\begin{defn}\label{primitive}
Let \ $\Lambda$ \ be a subset of \ $\mathbb{Q} \, \cap \, ]0,1]$. We say that a regular (a,b)-module \ $E$ \ is {\bf \ $[\Lambda]-$primitive} when all roots of its Bernstein polynomial are in \ $- \Lambda + \mathbb{Z}$.
\end{defn}

The following easy proposition is proved in [B.09]

\begin{prop}\label{primitive part}
Let \ $E$ \ be a regular (a,b)-module and \ $\Lambda$ \ a subset of \ $\mathbb{Q} \, \cap \, ]0,1]$. Then there exists a maximal submodule \ $E_{[\Lambda]}$ \ in \ $E$ \ which is \ $[\Lambda]-$primitive. Moreover the quotient \ $E\big/E_{[\Lambda]}$ \ is a \ $[\Lambda^c]-$primitive (a,b)-module, where we denote \  $\Lambda^c : = \mathbb{Q} \, \cap \, ]0,1] \setminus \Lambda$.
\end{prop}

We shall mainly consider the case where \ $\Lambda$ \ is a single element. Note that  the \ $[\lambda]-$primitive part of an (a,b)-module  \ $E \subset \Xi^{(N)}_{\Lambda}\otimes V$ \ corresponds to its intersection  with \ $\Xi^{(N)}_{\lambda}\otimes V$.

\subsection{The principal Jordan-H{\"o}lder sequence.}

The classification of rank \ $1$ \ regular (a,b)-modules is very simple : each isomorphy class is given by a complex number and to \ $\lambda \in \C$ \ corresponds the isomorphy class of \ $E_{\lambda} : = \A\big/\A.(a - \lambda.b)$. Then a {\bf Jordan-H{\"o}lder sequence} for a regular (a,b)-module \ $E$ \ is a sequence 
$$ 0 = F_0 \subset F_1 \subset \dots \subset F_k = E  $$
of normal sub-modules such that for each \ $j \in [1,k]$ \ the quotient \ $F_j\big/F_{j-1}$ \ has rank \ $1$. Then to each J-H. sequence we may associate an ordered sequence of  complex numbers \ $\lambda_1, \dots, \lambda_k$ \ such that \ $F_j\big/F_{j-1} \simeq E_{\lambda_j}$.

\parag{Example}
A regular rank 1 (a,b)-module is a fresco if and only if it is isomorphic to \ $E_{\lambda}$ \ for some \ $\lambda \in \mathbb{Q}^{+*}$. All rank 1 frescos are themes.
The classification of rank 2 regular (a,b)-modules obtained in [B.93] gives the list of \ $[\lambda]-$primitive rank 2 frescos  which is the following, where \ $\lambda_1 > 1$ \ is a rational number : 
\begin{align*}
& E =  E \simeq \A\big/\A.(a -\lambda_1.b).(a - (\lambda_1-1).b)  \tag{1}\\
& E \simeq \A\big/\A.(a -\lambda_1.b).(1+ \alpha.b^p)^{-1}.(a - (\lambda_1 + p-1).b) \tag{2}
\end{align*}
where \ $p \in \mathbb{N} \setminus \{0\}$ \ and \ $\alpha \in \C$.\\
The themes in this list are these in \ $(1)$ \ and these in \ $(2)$ \ with \ $\alpha \not= 0 $. For a \ $[\lambda]-$primitive fresco  in case \ $(2)$ \ the number \ $\alpha$ \ will be called the {\bf parameter} of the rank \ $2$ \ fresco. By convention we shall define \ $\alpha : = 1$ \ in the case \ $(1)$.$\hfill \square$ \\

The following existence result is proved in [B.93]

\begin{prop}\label{J-H exist}
For any regular (a,b)-module \ $E$ \ of rank \ $k$ \ there exists a J-H. sequence. The numbers \ $exp(2i\pi.\lambda_j)$ \ are independent of the J-H. sequence, up to permutation. Moreover the number \ $\mu(E) : = \sum_{j=1}^k \ \lambda_j $ \ is also independent of the choice of the J-H. sequence of \ $E$.
\end{prop}

\parag{Exercice} Let \ $E$ \ be a regular (a,b)-module and  \ $E' \subset E$ \ be a sub-(a,b)-module with the same rank than \ $E$. Show that \ $E'$ \ has finite \ $\C-$codimension in \ $E$ \ given by
$$ \dim_{\C}E\big/E' = \mu(E') - \mu(E) .$$
{\em hint} : make an induction on the rank of \ $E$.$\hfill \square$\\ 

For a \ $[\lambda]-$primitive fresco a more precise result is proved in [B.09].
 The following  generalization to the case of any fresco  is an easy application of the proposition \ref{primitive part}.
 
\begin{prop}\label{principal}
For any J-H. sequence of a rank \ $k$  \  fresco \ $E$ \  the numbers \ $\lambda_j+j-k, j \in [1,k]$ \ are the opposite of the  roots (with multiplicities) of the Bernstein polynomial of \ $E$. So, up to a permutation, they are independent of the choice of the J-H. sequence. Moreover, when \ $E$ \ is \ $[\lambda]-$primitive, there always exists a J-H. sequence such that the associated  sequence \ $\lambda_j+j$ \ is increasing.
\end{prop}

For a \ $[\lambda]-$primitive theme the situation is extremely rigid (see [B.09]) : 

\begin{prop}\label{theme J-H.}
Let \ $E$ \ a rank \ $k$ \ $[\lambda]-$primitive theme. Then, for each \ $j$ \ in \  $ [0,k]$,  there exists an {\bf unique} normal rank \ $j$ \ submodule \ $F_j$. The numbers associated to the unique J-H. sequence satisfy the condition that \ $\lambda_j+j, j \in [1,k]$ \ is a increasing  sequence.
\end{prop}

\parag{Order on \ $\mathbb{Q}$} We shall fix a total order on \ $\mathbb{Q}\big/\mathbb{Z}$, for instance the order induced by the inclusion \ $\mathbb{Q}\big/\mathbb{Z} \subset ]0,1]$. But any other total order would be also convenient. Then we shall use the corresponding order on \ $\mathbb{Q}$ \ defined as follows :\\
\begin{enumerate} [i)]
\item If \ $\lambda, \mu $ \  are in \ $\mathbb{Q}$ \ such \ $[\lambda] < [\mu]$ \ in \ $\mathbb{Q}\big/\mathbb{Z}$ \ we put \ $\lambda \tilde{<} \mu$.
\item If we have \ $[\lambda] = [\mu] $ \ in \ $\mathbb{Q}$ \ we put \ $\lambda \tilde{<} \mu$ \ when we have \ $\lambda < \mu$ \ for the usual order of  \ $\mathbb{Q}$.
\end{enumerate}
Remark that with such an order, a increasing sequence \ $\lambda_1, \dots, \lambda_k$ \ is such that theses numbers appears successively in the different classes modulo \ $\mathbb{Z}$ \ in increasing order in \ $\mathbb{Q}\big/\mathbb{Z}$ \ and in each class they are increasing in the usual sens. \\
Note that using the proposition \ref{primitive part} it is easy to see that any rank \ $k$ \  fresco admits a J.H. sequence such corresponding sequence of rational numbers \ $\lambda_1, \dots, \lambda_k$ satisfies the condition that the sequence \ $[\lambda_1], \dots, [\lambda_k]$ \ is increasing in \ $\mathbb{Q}\big/\mathbb{Z}$.

\begin{defn}\label{Principal J-H.}
Let \ $E$ \ be a  fresco of rank \ $k$ \ and let  
 $$0 = F_0 \subset F_1 \subset \dots \subset F_k = E$$
 be a J-H. sequence of \ $E$. Then for each \ $j \in [1,k]$ \ we have \ $F_j\big/F_{j-1} \simeq E_{\lambda_j}$, where \ $\lambda_1, \dots, \lambda_k$ \ are in \ $\mathbb{Q}^{+*}$.
We shall say that such a J-H. sequence is {\bf principal} when the sequence \  $[1,k]\ni j \mapsto \lambda_j + j$ \ is increasing for the order of \ $\mathbb{Q}$ \ fixed above.
\end{defn}

\begin{thm}\label{Uniqueness}
Let \ $E$ \ be a  fresco. Then \ $E$ \ admits an {\bf unique} principal J-H. sequence.
\end{thm}

 \parag{proof} The existence is an easy consequence of the propositions \ref{primitive part} and \ref{principal}. As we know the uniqueness of the \ $[\lambda]-$primitive part for any regular (a,b)-module, and that for a fresco \ $E$ \ the \ $[\lambda]-$primitive part is again a fresco (and normal so that the quotient is again a fresco), it is enough to prove the uniqueness in the \ $[\lambda]-$primitive case. We shall prove the uniqueness by induction on the rank \ $k$ \ of the \ $[\lambda]-$primitive fresco  \ $E$. \\

We begin by  the case of rank 2.

\begin{lemma}\label{petit}
Let \ $E$ \ be a rank 2 \ $[\lambda]-$primitive fresco and let \ $\lambda_1,\lambda_2$ \ the numbers corresponding to a principal J-H. sequence of \ $E$ \ (so  \ $\lambda_1+1 \leq \lambda_2+2$). Then the normal rank 1 submodule of \ $E$ \ isomorphic to \ $E_{\lambda_1}$ \ is unique. Moreover, if there exists a \ $\mu \not= \lambda_{1}$ \ and a rank \ $1$ \ normal submodule isomorphic to \ $E_{\mu}$ \ then \ $\mu = \lambda_{2}+ 1$. In this case there exists infinitely many different normal rank \ $1$ \ submodules isomorphic to \ $E_{\lambda_{2}+1}$.
\end{lemma}

\parag{Proof} The case \ $\lambda_1+ 1 = \lambda_2+2$ \ is obvious because then \ $E$ \ is a \ $[\lambda]-$primitive theme (see the example at the beginning of the paragraph 1.2 or  [B.10] corollary 2.1.7). So we may assume that \ $\lambda_2 = \lambda_1 + p_1-1$ \ with \ $p_1 \geq 1$ \ and that  \ $E$ \ is the quotient \ $E \simeq \A\big/\A.(a - \lambda_1.b).(a - \lambda_2.b)$ \ (see the classification of rank 2 frescos in 2.2), because the result is clear when \ $E$ \ is a theme.
We shall use the \ $\C[[b]]-$basis \ $e_1, e_2$ \  of \ $E$ \ where \ $a$ \ is defined by the relations
$$ (a - \lambda_2.b).e_2 = e_1 \quad {\rm and} \quad  (a - \lambda_1.b).e_1 = 0. $$
This basis comes from the isomorphism \ $E \simeq \A\big/\A.(a - \lambda_1.b).(a - \lambda_2.b)$ \  deduced from the classification of rank 2 frescos with \ $e_2 = [1]$ \ and \ $e_1 = (a - \lambda_2.b).e_2$.\\
Let look for \ $x : = U.e_2 + V.e_1$ \ such that \ $(a - \mu.b).x = 0$. Then we obtain
$$ b^2.U'.e_2 + U.(a - \lambda_2.b).e_2 + (\lambda_2- \mu).b.U.e_2 + b^2.V'.e_1 + (\lambda_{1} - \mu).b.V.e_{1}= 0$$
which is equivalent to the two equations :
$$ b^2.U' + (\lambda_2- \mu).b.U = 0 \quad {\rm and} \quad U + b^2.V' + (\lambda_{1} - \mu).b.V = 0 .$$
The first equation gives \ $U = 0$ \ for \ $\mu \not\in \lambda_{2} + \mathbb{N}$. As the case \ $U = 0$ \ will give (as we want also that \ $x \not\in b.E$) that \ $ x $ \ is equal to \ $ e_{1}$, up to a non zero multiplicative constant, we may assume that \ $\mu = \lambda_{2} + q$ \ for some \ $ q \in \mathbb{N}$. Moreover,  as the second equation implies \ $U(0) = 0$, we may assume that \ $q \geq 1$. This already shows that \ $ \mu \not= \lambda_{1}$ \ ($\lambda_{2} +1 = \lambda_{1}+  p_{1} > \lambda_{1}$) and this proves the first assertion. Now to finish our computation of normal rank \ $1$ \ sub-modules, we have \ $U  = \rho.b^{q}$. Then the solutions in \ $\C[[b]]$ \ of the equation
$$ \rho.b^{q-1} + b.V' - (p_{1}+q-1).V = 0  $$
are given by : $ V = -(\rho/p).b^{q-1} + \sigma.b^{p_{1}+q-1} $ \ and the condition  \ $x \not\in b.E$ \ implies now \ $ q =1$. So we obtain \ $\mu = \lambda_{2}+1$ \ and for each \ $\tau \in \C$ \ the element \ $x = (1 - p_{1}.\tau.b^{p_{1}}).e_{1} - p_{1}.b.e_{2}$ \ generates a normal rank \ $1$ \ sub-module isomorphic to \ $E_{\lambda_{2}+1}$. And with the unique sub-module isomorphic to \ $E_{\lambda_{1}}$ \ they are all the normal rank \ $1$ \ sub-modules in such an \ $E$. $\hfill \blacksquare$\\

\parag{end of the proof of theorem \ref{Uniqueness}} As the result is obvious for \ $k = 1$, we may assume \ $k \geq 2$ \ and the result proved in rank \ $\leq k-1$. Let \ $F_j, j \in [1,k]$ \ and \ $G_j, j \in [1,k]$ \ two J-H. principal sequences for \ $E$. As the sequences \ $\lambda_j + j$ \ and \ $\mu_j+j$ \ co{\"i}ncide up to the order (they are of the form \ $ -x_{j}  + k$ \ where the \ $(x_{j})_{j \in[1,k]}$ \ are the roots of the Bernstein polynomial, counting multiplicities) and are both increasing, they  co{\"i}ncide. Now let \ $j_0$ \ be the first integer in \ $[1,k]$ \ such that \ $F_{j_0} \not= G_{j_0}$. If \ $j_0 \geq 2$ \ applying the induction hypothesis to \ $E\big/F_{j_0-1}$ \ gives \ $F_{j_0}\big/F_{j_0-1} = G_{j_0}\big/F_{j_0-1}$ \ and so \ $F_{j_0} = G_{j_0}$. \\
So we may assume that \ $j_0 = 1$. Let \ $H$ \ be the normalization of \ $F_1 + G_1$. As \ $F_1$ \ and \ $G_1$ \ are normal rank 1 and distinct, then \ $H$ \ is a rank 2 normal sub-module. It is a \ $[\lambda]-$primitive fresco of rank 2 with two different normal rank 1 sub-modules which are isomorphic as \ $\lambda_1 = \mu_1$. Moreover the principal J-H. sequence of \ $H$ \ begins by a normal submodule isomorphic to \ $E_{\lambda_1}$. This contradicts the previous lemma and so  \ $F_1 = G_1$. So for any \ $j \in [1,k]$ \ we have \ $F_j = G_j$. $\hfill \blacksquare$

\parag{Notation} For a fresco \ $E$ \ we shall denote \ $F_{j}(E)$ \ the \ $j-$th term of the principal J.H. sequence of \ $E$. Remark that the isomorphism class of the fresco \ $F_{j}(E)$ \ depends only on the isomorphism class of \ $E$, thanks to the previous theorem : it implies that any isomorphism \ $f : E \to E'$ \  of frescos sends \  $F_{j}(E)$ \ on  \ $F_{j}(E')$.\\
Remark also that for any change of variable\footnote{see [B.11].} \ $\theta \in \C[[a]]$ \ we have a natural isomorphism  \ $\theta_{*}(F_{j}(E)) \simeq F_{j}(\theta_{*}(E))$, again as a consequence of the previous theorem and the invariance of the Bernstein polynomial by changes of variable.\\

\begin{defn}\label{Inv. fond.}
Let \ $E$ \ be a  fresco and consider its principal J-H. sequence \ $F_j, j \in [1.k]$. Put \ $F_j\big/F_{j-1} \simeq E_{\lambda_j}$ \ for \ $j \in [1,k]$ \ (with \ $F_0 = \{0\}$). We shall call the {\bf fundamental invariants} of \ $E$ \ the ordered k-tuple\ $(\lambda_1, \dots, \lambda_k) $.
\end{defn}

Of course, if we have any J-H. sequence for a  fresco  \ $E$ \ with rank \ $1$ \ quotients associated to the rational numbers \ $\mu_1, \dots, \mu_k$, it is easy to recover the fundamental invariants of \ $E$ \ because the numbers \ $\mu_j+j, j \in [1,k]$ \ are the same than the numbers \ $\lambda_j+j$ \ up to a permutation. But as the sequence \ $\lambda_j+j$ \ is increasing  with \ $j$ \ for the fixed order on \ $\mathbb{Q}$, it is enough to put the \ $\mu_j+j$ \ in the increasing order with \ $j$ \  to recover the fundamental invariants of \ $E$.

\begin{cor}\label{gen. Fj}
Let \ $\lambda_1, \dots, \lambda_k$ \ the fundamental invariants of a rank \ $k$ \  fresco \ $E$ \ and assume that we have an isomorphism
$$ E \simeq \A\big/\A.P \quad {\rm where} \quad  P : = (a - \lambda_1.b).S_1^{-1} \dots S_{k-1}^{-1}.(a - \lambda_k.b).S_k^{-1} $$
where \ $S_1, \dots, S_{k}$ \ are invertible elements in \ $\C[[b]]$,  by sending the generator \ $e \in E$ \ to \ $1 \in \A$. Then for each \ $j \in [1,k-1]$ \ the element
$$ e_j : = S_j^{-1}.(a - \lambda_{j+1}.b) \dots S_{k-1}^{-1}.(a- \lambda_k.b).S_k^{-1}.e $$
is a generator of the fresco \ $F_j$ \ which is the \ $j-$th term in the principal J.H. sequence of \ $E$.
\end{cor}

\parag{proof} By an easy induction, it is enough to prove that \ $S_{k-1}.e_{k-1} = (a - \lambda_{k}.b).e_{k}$ \ is a generator of the fresco  \ $F_{k-1}$. So define \ $G : = \A.(a - \lambda_k.b).e_{k}$. Let first prove that \ $G$ \ is normal in \ $E$. So assume that for some \ $x \in E$ \ we have \ $b.x \in G$. We may write
$$  x = u.e_{k} \quad {\rm and} \quad b.x = v.(a - \lambda_k.b).e_{k} $$
for some \ $u,v \in \A$. So we have \ $b.u - v.(a - \lambda_k.b) \in \A.P.S_{k} $. This implies that \ $b.u$ \ is in \ $\A.(a - \lambda_k.b)$. But \ $\A\big/\A.(a - \lambda_k.b) \simeq E_{\lambda_k}$ \ has no $b-$torsion. Then we have \ $u = w.(a - \lambda_k.b)$ \ and \ $x = w.(a - \lambda_k.b).e_{k}$. So \ $G$ \ is normal and \ $E\big/G$ \ is isomorphic to \ $E_{\lambda_k}$. Then the uniqueness of the principal J.H. sequence of \ $E$ \ implies the equality \ $G = F_{k-1}$. $\hfill \blacksquare$\\

So in the situation of the previous corollary we have a natural \ $\C[[b]]-$ basis \ $e_1, \dots ,e_k$ \ of \ $E$ \ associated to a generator  \ $e$ \ of \ $E$. Note that the uniqueness of \ $F_{k-1}$ \ implies the uniqueness of \ $S_k$ \ up to a non zero constant and then we have the uniqueness of \ $S_{k-1}$ \ again up to a non zero constant  etc$\dots$. So \ $P$ \ is unique up to  a non zero constant also (and we may normalize it asking that the coefficient of \ $a^k$ \ is \ $1$). This means that the elements \ $e_{1}, \dots, e_k$ \ are well determined by the choice of \ $e$ \ up to a diagonal invertible matrix with coefficients in \ $\C$. And for each \ $j \in [1,k]$ \ we have \ $F_j = \A.e_j$.\\

\bigskip

Let \ $\lambda_1, \dots, \lambda_k$ \ be the fundamental invariants of a rank \ $k$ \  fresco \ $E_0$, and note \ $\mathcal{F}(\lambda_1, \dots, \lambda_k)$ \ the set of isomorphism classes of frescos with these fundamental invariants. The uniqueness of the principal J-H. sequence of a fresco allows to define for each \ $(i,j) \ 1 \leq i < j \leq k$ \ a map
$$ q_{i,j} : \mathcal{F}(\lambda_1, \dots, \lambda_k) \to \mathcal{F}(\lambda_i, \dots, \lambda_j)$$
defined by \ $q_{i,j}([E]) = [F_j\big/F_{i-1}] $ \ where \ $(F_h)_{h\in [0,k]}$ \ is the principal J-H. sequence of \ $E$. This makes sense because we know that  any isomorphism \ $\varphi : E^1 \to E^2$ \ between two  frescos induces isomorphisms between each term of  the corresponding principal J-H. sequences.\\
Now from the corollary above it is clear that such an application is algebraic (see the definition in section 3).\\

For instance  the classification of rank \ $2$ \ $[\lambda]-$primitive frescos gives (see example before proposition \ref{J-H exist}) for any rational number \ $\lambda_1 > 1$ \ and \ $p_1 \in \mathbb{N}$ \ the following description :
\begin{itemize}
\item ${\rm for} \quad  p_1 = 0 \quad  \mathcal{F}(\lambda_1, \lambda_1-1) \simeq \{pt\} $.
\item ${\rm for} \quad p_1 \geq 1 \quad \alpha : \mathcal{F}(\lambda_1, \lambda_1 + p_1 - 1) \overset{\simeq}{\longrightarrow}  \C .$
\end{itemize}
where \ $\{pt\}$ \ is given by the isomorphism class of \ $\A\big/\A.(a - \lambda_1.b).(a - (\lambda_1-1).b) $ \ and where the isomorphism class associated to \ $\alpha^{-1}(x)$ \ for \ $x  \in \C$ \ in the case \ $p_1 \geq 1$ \ is given by \ $\A\big/\A.(a - \lambda_1.b).(1 + x.b^{p_1})^{-1}.(a - (\lambda_1 + p_1 -1).b)$. \\

This implies that the function \ $\alpha_j : \mathcal{F}(\lambda_1, \dots, \lambda_k) \to \C$ \ given by sending \ $[E]$ \ to the parameter \ $\alpha([F_{j+1}\big/F_{j-1}])$ \ for \ $j \in [1, k-1]$ \ is a polynomial (see section 3 for the definition). Remark that this function is not identically \ $0$ \ or \ $1$ \ if and only if \ $\lambda_{j+1} - \lambda_j+ 1$ \ is a positive integer.

\section{Semi-simplicity.}

\subsection{Semi-simple regular  (a,b)-modules.}

\begin{defn}\label{ss 0}
Let \ $E$ \ be a regular (a,b)-module. We say that \ $E$ \ is {\bf semi-simple} if it is a sub-module of a finite direct sum of rank 1 regular (a,b)-modules.
\end{defn}

Note that if \ $E$ \ is a sub-module of a regular (a,b)-module it is necessary regular. As a direct sum of regular (a,b)-modules is regular, our assumption that \ $E$ \ is regular is superfluous.\\
It is clear from this definition that a sub-module of a semi-simple (a,b)-module is semi-simple and that a (finite) direct sum of semi-simple (a,b)-modules is again semi-simple.

\parag{Remark} A rather easy consequence of the classification of rank \ $2$ \ (a,b)-modules given in [B.93]  is that the rank \ $2$ \ simple pole (a,b)-modules defined in the \ $\C[[b]]-$basis \ $x,y$ \  by the relations :
$$(a - (\lambda+p).b).x = b^{p+1}.y \quad {\rm and} \quad (a - \lambda).y = 0$$
for any \ $\lambda \in \C$ \ and \ any \ $p \in \mathbb{N}$ \ are not semi-simple. We leave the verification of this point to the reader.$\hfill\square$\\

Let us begin by a characterization of the semi-simple (a,b)-modules which have a simple pole. First we shall prove that a quotient of a semi-simple (a,b)-module is semi-simple. This will be deduced from the following lemma.

\begin{lemma}\label{quot. ss 1}
Let \ $E$ \ be a (a,b)-module which is direct sum of regular rank \ $1$ \ (a,b)-modules, and let \ $F \subset E$ \ be a rank  $1$ \ normal sub-module. Then \ $F$ \ is a direct factor of \ $E$.
\end{lemma}

\begin{cor}\label{quot. ss 2}
If \ $E$ \ is a semi-simple regular (a,b)-module and \ $F$ \ a normal sub-module of \ $E$, the quotient \ $E\big/F$ \ is a (regular) semi-simple (a,b)-module.
\end{cor}

\parag{Proof of the lemma} Let  \ $E = \oplus_{j=1}^k \ E_{\lambda_j}$ \ and assume that \ $F \simeq E_{\mu}$. Let \ $e_j$ \ be a standard generator of \ $E_{\lambda_j}$ \ and \ $e$ \ a standard generator of \ $E_{\mu}$. Write
$$ e = \sum_{j=1}^k \ S_j(b).e_j $$
and compute \ $(a - \mu.b).e = 0$ \ using the fact that \ $e_j, j \in [1,k]$ \ is a \ $\C[[b]]-$basis of \ $E$ \ and the relations \ $(a - \lambda_j.b).e_j = 0$ \ for each \ $j$. We obtain for each \ $j \in [1,k]$ \ the relation
$$ b.S'_j - (\mu-\lambda_j).S_j = 0 .$$
So, if \ $\mu - \lambda_j$ \ is not in \ $\mathbb{N}$, we have \ $S_j = 0$. When \ $\mu = \lambda_j + p_j$ \ with \ $p_j \in \mathbb{N}$ \ we obtain \ $S_j(b) = \rho_j.b^{p_j}$ \ for some \ $\rho_j \in \C$. As we assume that \ $e$ \ is not in \ $b.E$, there exists at least one \ $j_0 \in [1,k]$ \ such that \ $p_{j_0} = 0$ \ and \ $\rho_{j_0} \not= 0 $. Then it is clear that we have
$$ E = F \oplus \big(\oplus_{j\not= j_0} E_{\lambda_j}\big) $$
concluding the proof. $\hfill \blacksquare$

\parag{Proof of the corollary} We argue by induction on the rank of \ $F$. In the rank \ $1$ \ case, we have \ $F \subset E \subset \mathcal{E} : = \oplus_{j=1}^k \ E_{\lambda_j}$. Let \ $\tilde{F}$ \ the normalization of \ $F$ \ in \ $\oplus_{j=1}^k \ E_{\lambda_j}$. Then the lemma shows that there exists a \ $j_0 \in [1,k]$ \ such that 
 $$\mathcal{E} = \tilde{F} \oplus\big( \oplus_{j\not= j_0} E_{\lambda_j}\big).$$
 Then, as \ $\tilde{F} \cap E = F$, the quotient map \ $\mathcal{E} \to \mathcal{E}\big/\tilde{F} \simeq \oplus_{j\not= j_0} E_{\lambda_j}$ \ induces an injection of \ $E\big/F$ \ in a direct sum of regular rank \ $1$ \ (a,b)-modules. So \ $E\big/F$ \ is semi-simple.\\
 Assume now that the result is proved for \ $F$ \ with rank \ $\leq d-1$ \ and assume that \ $F$ \ has rank \ $d$. Then using a rank \ $1$ \ normal sub-module \ $G$ \ in \ $F$, we obtain that \ $F\big/G $ \ is a normal rank \ $d-1$ \ sub-module of \ $E\big/G$. Using the rank \ $1$ \ case we know that \ $E\big/G$ \ is semi-simple, and the induction hypothesis gives that 
  $$E\big/F = (E\big/G)\Big/(F\big/G)$$
   is semi-simple. $\hfill \blacksquare$\\

\begin{prop}\label{simple pole ss}
Let \ $E$ \ be a simple pole semi-simple (a,b)-module. Then \ $E$ \ is  a direct sum of regular rank \ $1$ \ (a,b)-modules.
\end{prop}

\parag{proof}  We shall prove the proposition by induction on the rank \ $k$ \ of \ $E$. As the rank \ $1$ \ case is obvious, assume the proposition proved for \ $k-1$ \ and  that \ $E$ \ is a simple pole rank \ $k \geq 2$ \ semi-simple (a,b)-module. From the existence of J-H. sequence, we may find an exact sequence
$$ 0  \to F \to E \to  E_{\lambda} \to  0 $$
where \ $F$ \ has rank \ $k-1$. By definition \ $F$ \ is semi-simple, but it also has a simple pole because \ $a.F \subset a.E \subset b.E $ \ implies that \ $a.F \subset F \cap b.E = b.F$ \ as \ $F$ \ is normal in \ $E$.
So by the induction hypothesis \ $F$ \ is a direct sum of  regular rank \ $1$ \ (a,b)-modules.\\
Let \ $e \in E$ \ such that its image in \ $E_{\lambda}$ \ is \ $e_{\lambda}$ \ a standard generator of \ $E_{\lambda}$ \ satisfying \ $a, e_{\lambda} = \lambda.b.e_{\lambda}$. Then we have \ $(a - \lambda.b).e \in F$. 
\smallskip

We shall first look at the case \ $k = 2$. 

\smallskip

So \ $F$ \ is a rank \ $1$ \ and  we have \ $F \simeq E_{\mu}$ \ for some \ $\mu \in \C$. Let \ $e_{\mu}$ \ be a standard generator in \ $F$ \ and put
$$ (a - \lambda.b).e = S(b).e_{\mu}.$$
 Our simple pole assumption on \ $E$ \ implies \ $S(0) = 0$ \ and we may write \\
  $S(b) = b.\tilde{S}(b)$. We look for \ $T \in \C[[b]]$ \ such that \ $\varepsilon : = e + T(b).e_{\mu}$ \ satisfies
$$ (a - \lambda.b).\varepsilon = 0.$$
If such a \ $T$ \ exists, then we would have \ $E = E_{\mu}\oplus E_{\lambda}$ \ where \ $E_{\lambda} $ \ is the sub-module generated by \ $\varepsilon$, because it is clear that we have \ $E = F \oplus \C[[b]].e$ \ as a \ $\C[[b]]-$module. \\
To find \ $T$ \ we have to solve in \ $\C[[b]]$ \  the differential equation
$$ b.\tilde{S}(b) + (\mu - \lambda).b.T(b) + b^2.T'(b) = 0 .$$
If \ $\lambda- \mu$ \ is not a non negative integer, such a \ $T$ \ exists and is unique. But when \ $\lambda = \mu + p$ \ with \ $p \in \mathbb{N}$, the solution exists if and only if the coefficient of \ $b^p$ \ in \ $\tilde{S}$ \ is zero. If it is not the case, define \ $\tilde{T}$ \ as the solution of the differential equation
$$ \tilde{S}(b) -  \alpha.b^p + b.\tilde{T}'(b) - p.\tilde{T}(b) = 0 $$
where \ $\alpha \not= 0$ \ is the coefficient of \ $b^p$ \ in \ $\tilde{S}$ \ and where we choose \ $\tilde{T}$ \ by asking that it has no \ $b^p$ \ term. Then \ $\varepsilon_1 : = e + \tilde{T}(b).e_{\lambda-p}$ \ satisfies
$$ (a - (\lambda - p).b).\varepsilon_1 = \alpha.b^{p+1}.e_{\lambda -p} .$$
Then, after changing \ $\lambda-p$ \ in \ $\lambda$ \ and  \ $e_{\lambda-p}$ \ in \ $\alpha.e_{\lambda-p}$,  we recognize one of the  rank \ $2$ \ modules  which appears in the remark following the definition \ref{ss 0},  and which is not semi-simple. So we have a contradiction. This concludes the rank \ $2$ \ case.

\bigskip

Now consider the case \ $k \geq  3$ \ and using the induction hypothesis write
$$ F \simeq \oplus_{j=1}^{k-1} \ E_{\mu_j} $$
and denote \ $e_j$ \ a standard generator for \ $E_{\mu_j}$. Write
$$ (a - \lambda.b).e = \sum_{j=1}^{k-1} \ S_j(b).e_j  .$$
The simple pole assumption again gives \ $S_j(0) = 0 $ \ for each \ $j$; now we 
look for \ $T_1, \dots ,T_{k-1}$ \ in \ $\C[[b]]$ \ such that, defining \ $\varepsilon : = e +  \sum_{j=1}^{k-1} \ T_j(b).e_j $, we have
$$ (a - \lambda.b).\varepsilon =  0  .$$
For each \ $j$ \ this leads to the differential equation
$$b.\tilde{S}_j(b) + (\mu_j - \lambda).b.T_j(b) + b^2.T_j'(b) = 0 $$
where we put \ $S_j(b) = b.\tilde{S}_j(b)$. We are back, for each given \ $j \in [1,k-1]$, to the previous problem in the rank \ $2$ \ case. So if for each \ $j \in [1,k-1]$ \ we have a solution \ $T_j \in \C[[b]]$ \  it is easy to conclude that
$$ E \subset F \oplus \C[[b]].\varepsilon \simeq F \oplus E_{\lambda}.$$
If for some \ $j_0$ \ there is no solution, the coefficient of \ $b^{j_0}$ \ in \ $\tilde{S}_{j_0}$ \ does not vanish  and then the image of \ $E$ \ in \ $E\big/\oplus_{j \not= j_0} E_{\mu_j} $ \ is a rank \ $2$ \ not semi-simple (a,b)-module. This contradicts the corollary \ref{quot. ss 2} and concludes the proof. $\hfill \blacksquare$\\

\parag{Remark} Let \ $E$ \ be a semi-simple  rank \ $k$ \ (a,b)-module and let \ $E^{\sharp}$ \ its saturation by \ $b^{-1}.a$, then \ $E^{\sharp}$ \ is semi-simple because there exists \ $N \in \mathbb{N}$ \  with an inclusion \  $E^{\sharp} \subset b^{-N}.E $. Also \ $E^{\flat}$, the maximal simple pole sub-module of \ $E$, is semi-simple. Then we have 
$$ E^{\flat} \simeq \oplus_{j=1}^k \ E_{\lambda_j+d_j} \subset E \subset E^{\sharp} \simeq \oplus_{j=1}^k \ E_{\lambda_j} $$
where \ $d_1, \dots, d_k$ \ are non negative integers.\\

\begin{cor}\label{dual ss ; quot. ss}
Let \ $E$ \ be a semi-simple (a,b)-module. Then its dual \ $E^*$ \ is semi-simple. 
\end{cor}

\parag{proof} For a regular (a,b)-module \ $E$ \ we have \ $(E^{\flat})^* = (E^*)^{\sharp}$. But the dual of \ $E_{\lambda}$ \ is \ $E_{-\lambda}$ \ and the duality commutes with direct sums. So we conclude that \ $E^* \subset (E^{\flat})^*$ \ and so \ $E^*$ is a semi-simple. $\hfill \blacksquare$

\parag{Remark} The tensor product of two semi-simple (a,b)-modules is semi-simple as a consequence of the fact that \ $E_{\lambda} \otimes E_{\mu} \simeq E_{\lambda+\mu}$.

\subsection{The semi-simple filtration.}

\bigskip

\begin{defn}\label{element ss 1}
Let \ $E$ \ be a regular (a,b)-module and \ $x$ \ an element in \ $E$. We shall say that \ $x$ \ is {\bf semi-simple} if \ $\A.x$ \ is a semi-simple (a,b)-module.
\end{defn}

It is clear that any element in a semi-simple  (a,b)-module is semi-simple. The next lemma shows that the converse is true.

\begin{lemma}\label{element ss 2}
Let \ $E$ \ be a regular (a,b)-module such that any \ $x \in E$ \ is semi-simple. Then \ $E$ \ is semi-simple.
\end{lemma}

\parag{proof} Let \ $e_1, \dots, e_k$ \ be a \ $\C[[b]]-$basis of \ $E$. Then each \ $\A.e_j$ \ is semi-simple, and the map \ $\oplus_{j=1}^k \A.e_j \to E $ \ is surjective. So \ $E$ \ is semi-simple thanks to the corollary \ref{quot. ss 2} and the comment following the definition \ref{ss 0}. $\hfill \blacksquare$\\

\begin{lemma}\label{ss part 0}
Let \ $E$ \ be a regular (a,b)-module. The subset \ $S_1(E)$ \ of semi-simple elements in \ $E$ \ is a normal submodule in \ $E$.
\end{lemma}

\parag{proof} As a direct sum and quotient of semi-simple (a,b)-modules are semi-simple, it is clear that for  \ $x$ \ and \ $y$ \ semi-simple the sum \ $\A.x + \A.y$ \ is semi-simple. So \ $x + y$ \ is semi-simple. This implies that \ $S_1(E)$ \ is a sub-module of \ $E$. If \ $b.x$ \ is in \ $S_1(E)$, then \ $\A.b.x$ \ is semi-simple. Then \ $\A.x \simeq \A.b.x \otimes E_{-1}$ \ is also semi-simple, and so \ $S_1(E)$ \ is normal in \ $E$. $\hfill \blacksquare$

\begin{defn}\label{ss part 1}
Let \ $E$ \ be a regular (a,b)-module. Then the sub-module \ $S_1(E)$ \ of semi-simple elements in \ $E$ \ will be called the {\bf semi-simple part} of \ $E$.\\
Defining \ $S_j(E)$ \ as the pull-back on \ $E$ \ of the semi-simple part of \ $E\big/ S_{j-1}(E)$ \ for \ $j \geq 1$ \ with the initial condition \ $S_0(E) = \{0\}$, we define a sequence of normal sub-modules in \ $E$ \ such that \ $S_j(E)\big/S_{j-1}(E) = S_1(E\big/S_{j-1}(E))$. We shall call it the {\bf \em semi-simple filtration} of \ $E$. The smallest integer \ $d$ \ such we have \ $S_d(E) = E$ \ will be called the {\bf semi-simple depth} of \ $E$ \ and we shall denote it \ $d(E)$.\\
\end{defn}

\parag{Remarks}
\begin{enumerate}[i)]
\item As \ $S_1(E)$ \ is the maximal semi-simple sub-module of \ $E$ \ it contains any rank \ $1$ \ sub-module of \ $E$. So \ $S_1(E) = \{0\}$ \ happens if and only if \ $E = \{0\}$.
\item Then the ss-filtration of \ $E$ \ is strictly increasing for \ $0 \leq j \leq d(E)$. 
\item It is easy to see that for any submodule \ $F$ \ in \ $E$ \ we have \ $S_j(F) = S_j(E) \cap F$. So \ $S_j(E)$ \ is the subset of \ $x \in E$ \ such that \ $ d(x) : = d(\A.x) \leq j$ \ and \  $d(E)$ \ is the supremum of \ $d(x)$ \ for \ $x \in E$. $\hfill \square$
\end{enumerate}

\begin{lemma}\label{primitive ss}
Let \ $0 \to F \to E \to G \to 0$ \ be an exact sequence of regular (a,b)-modules such that \ $F$ \ and \ $G$ \ are semi-simple and such that for any roots \ $\lambda, \mu$ \ of the Bernstein polynomial of \ $F$ \ and \ $G$ \ respectively we have \ $\lambda - \mu \not\in \mathbb{Z}$. Then \ $E$ \ is semi-simple.
\end{lemma}

\parag{proof} By an easy induction we may assume that \ $F$ \ is rank \ $1$ \ and also that \ $G$ \ is rank \ $1$. Then the conclusion follows from the classification of the rank \ $2$ \ regular (a,b)-module. $\hfill \blacksquare$

 \parag{Application} A regular (a,b)-module \ $E$ \ is semi-simple if and only if for each \ $\lambda \in \C$ \ its primitive part \ $E[\lambda]$ \ is semi-simple.\\
   
The following easy facts are  left as  exercices for the reader.

\begin{enumerate}
\item  Let \ $0 \to F \to E \to G \to 0$ \ be a short exact sequence of regular (a,b)-modules. Then we have the inequalities
  $$ \sup \{ d(F), d(G)\} \leq d(E) \leq d(F) + d(G) .$$
  \item Let \ $E \subset \Xi_{\lambda}^{(N)} \otimes V$ \ be a sub-module where \ $V$ \ is  a finite dimensional complex space. Then for each \ $j \in [1, N+1]$ \ we have
  $$ S_j(E) = E \cap \left( \Xi_{\lambda}^{(j-1)} \otimes V \right) .$$
\item  In the same situation as above,  the minimal \ $N$ \ for such an embedding of a geometric \ $[\lambda]-$primitive \ $E$ \ is equal to \ $d(E)-1$.
\item In the same situation as above, the equality  \ $S_1(E) = E \cap \Xi_{\lambda}^{(0)} \otimes V$ \ implies \ $rank(S_1(E)) \leq \dim_{\C} V $. It is easy to see that any \ $[\lambda]-$primitive semi-simple geometric (a,b)-module  \ $F$ \ may be embedded in \ $\Xi_{\lambda}^{(0)} \otimes V$ \ for some \ $V$ \ of dimension \ $rank(F)$. But it is also easy to prove that when \ $F = S_1(E)$ \ with \ $E$ \ as above, such an embedding may be extended to an embedding of \ $E$ \ into \ $\Xi_{\lambda}^{(d(E)-1)} \otimes V$. $\hfill \square$
\end{enumerate}

\subsection{The co-semi-simple filtration.}

We shall discuss briefly the filtration deduced by duality from the semi-simple filtration of a regular (a,b)-module.

\begin{lemma}\label{co-ss}
For  any regular (a,b)-module \ $E$ \ there exists a maximal co-semi-simple normal sub-module, noted \ $\Sigma_1(E)$. That is to say  that any normal sub-module \ $F$ \ of \ $E$ \ such that \ $E\big/F$ \ is semi-simple  contains \ $\Sigma_1(E)$, and that the quotient \ $E\big/\Sigma_{1}(E)$ \ is semi-simple.
\end{lemma}

\parag{Remark} We have \ $\Sigma_{1}(E) = 0$ \ if and only if \ $E$ \ is semi-simple. Note also that for \ $E \not= 0$ \ the inclusion \ $\Sigma_{1}(E) \subset E$ \ is strict because \ $E$ \ always has a co-rank \ $1$ \ normal sub-module, and any rank \ $1$ \ regular (a,b)-module is semi-simple. $\hfill \square$\\

\parag{proof} We shall prove that we have \ $ \Sigma_{1}(E) = (E^{*}\big/S_{1}(E^{*}))^{*}$. Dualizing the exact sequence of (a,b)-modules
$$ 0 \to S_{1}(E^{*}) \to E^{*} \to E^{*}\big/S_{1}(E^{*}) \to 0 $$
we see that \ $ G : = (E^{*}\big/S_{1}(E^{*}))^{*}$ \ is a normal submodule of \ $E$ \ and that \ $E\big/G$ \ is semi-simple because its dual \ $S_{1}(E^{*})^{*}$ \ is semi-simple.\\
Consider now a normal sub-module \ $H$ \ in \ $E$ \ such that \ $E\big/H$ \ is semi-simple. Then we have an exact sequence
$$0 \to (E\big/H)^{*} \to E^{*} \to H^{*} \to 0 $$
where the sub-module \ $(E\big/H)^{*}$ \ is semi-simple. So it is contained in \ $S_{1}(E^{*})$. The surjective map \ $H^{*} \simeq E^{*}\Big/(E\big/H)^{*}  \to G^{*}$ \ implies that \ $G \subset H$ \ proving the minimality of \ $G$. $\hfill \blacksquare$

\begin{defn}\label{co-1}
Let  \ $E$ \ a regular (a,b)-module. Define inductively, for any \ $h \geq 2$, the normal sub-modules \ $\Sigma_{h}(E) : = \Sigma_{1}(\Sigma_{h-1}(E))$. We shall call this filtration of \ $E$ \ the {\bf co-semi-simple filtration} of \ $E$.
\end{defn}

\begin{cor}\label{co-2}
For  any regular (a,b)-module \ $E$ \ and any \ $h \in [1, d(E^{*})]$ \ we have a natural isomorphism
$$ \Sigma_{h}(E) \simeq (E^{*}\big/S_{h}(E^{*}))^{*} .$$
\end{cor}

\parag{proof} As the case \ $h = 1$ \ has been proved in the previous lemma, assume \ $h \geq 2$ \ and the isomorphism obtained for \ $h-1$. Note \ $F : = \Sigma_{h-1}(E)$ ; the case \ $h = 1$ \  applied to \ $F$ \ gives \ $ \Sigma_{1}(F) \simeq (F^{*}\big/S_{1}(F^{*}))^{*}$. The induction hypothesis  implies
$$F^{*} \simeq E^{*}\big/S_{h-1}(E^{*}) $$
so we have \ $ \Sigma_{1}(F)^{*} \simeq E^{*}\big/S_{h-1}(E^{*})\Big/S_{1}(E^{*}\big/S_{h-1}(E^{*})) \simeq E^{*}\big/S_{h}(E^{*})$ \ and we conclude by dualizing this isomorphism. $\hfill \blacksquare$\\

An immediate consequence of this corollary is that we have \ $\Sigma_{d-1}(E) \not= 0$ \ and \ $\Sigma_{d}(E) = 0$ \ for \ $d = d(E^{*})$.

\begin{lemma}\label{co-3}
Let  \ $E$ \ a regular (a,b)-module. For any \ $h \geq 1$ \ we have
$$ \Sigma_{h}(E) \subset S_{d-h}(E)$$
with \ $d : = d(E)$.
\end{lemma}

\parag{proof} We shall make an induction on \ $h \geq 1$. As \ $E\big/S_{d-1}(E)$ \ is semi-simple, we  have \ $\Sigma_{1}(E) \subset S_{d-1}(E)$, by definition of \ $\Sigma_{1}(E)$.\\
Assume the inclusion proved for \ $h-1 \geq 1$. Then as \ $S_{d-h+1}(E)\big/S_{d-h}(E) = S_{1}(E\big/S_{d-h}(E))$ \ is semi-simple, the quotient
$$\Sigma_{h-1}(E)\Big/\Sigma_{h-1}(E) \cap S_{d-h}(E) \hookrightarrow S_{d-h+1}(E)\big/S_{d-h}(E) $$
is also semi-simple, and this implies that
$$\Sigma_{h}(E) = \Sigma_{1}(\Sigma_{h-1}(E)) \subset \Sigma_{h-1}(E) \cap S_{d-h}(E)  \subset   S_{d-h}(E) $$
which gives the result. $\hfill \blacksquare$

\begin{cor}\label{co-4}
For any regular (a,b)-module we have \ $d(E^{*}) =  d(E)$.
\end{cor}

\parag{proof} As \ $S_{1}(E)$ \ is semi-simple, $\Sigma_{d(E)-1}(E)$ \ is semi-simple, and then \ $\Sigma_{d(E)}(E) = 0$. Using the remark following the corollary \ref{co-2}, we deduce that \ $ d(E) \geq d(E^{*}) $. By duality, we obtain the equality. $\hfill \blacksquare$

\subsection{The case of frescos.}

In the case of a fresco, we shall give a simple characterization of the semi-simple filtration. Let me begin by a simple lemma.

\begin{lemma}
 \begin{itemize}
   \item If \ $T \subset E$ \ is a \ $[\lambda]-$primitive theme in a fresco \ $E$, its normalization is also a \ $[\lambda]-$primitive theme (of same rank than \ $T$).
   \item If \ $S \subset E$ \ is a semi-simple fresco in a fresco \ $E$, its normalization is also a semi-simple fresco.
 \end{itemize}
 \end{lemma}
 
 \parag{proof} Note first that a normal sub-module of a fresco is a fresco. So the normalization \ $\tilde{T}$ \ will be a \ $[\lambda]-$primitive theme if we can show that it contains a unique rank \ $1$ \ normal submodule\footnote{see [B. 10].}. But as \ $T$ \ has an unique normal rank \ $1$ \ submodule \ $L$, it is clear that the only rank \ $1$ \ normal sub-module of \ $\tilde{T}$ \ is the normalization \ $\tilde{L}$ \ of \ $L$.\\
 If now \ $S \subset E$ \ is a semi-simple fresco, then if a \ $\A-$linear map \ $ \varphi : \tilde{S} \to \Xi_{\lambda}$ \ has rank \ $ \geq 2$, then its restriction to \ $S$ \ has the same rank as \ $S$ \ has finite codimension in its normalzation \ $\tilde{S}$ \ and this contradicts the semi-simplicity of \ $S$. $\hfill\blacksquare$.\\
 
 For a fresco \ $E$ \  we have the following characterization for the semi-simple (and co-semi-simple) filtration of  \ $E$.

 \begin{prop}\label{ss sequence}
  Let \ $E$ \ be a fresco. Then we have the following properties :
  \begin{enumerate}[i)]
  \item Any \ $[\lambda]-$primitive  sub-theme \ $T$ \ in \ $E$ \ of rank \ $j$ \ is contained in \ $S_j(E)$.
  \item Any  \ $[\lambda]-$primitive quotient theme \ $T$ \ of \ $S_j(E)$ \ has  rank \ $\leq j$. 
  \item For any \ $j \in \mathbb{N}$ \ we have 
  $$ S_j(E) = \cap_{\lambda}\cap_{\varphi \in Hom_{\A}(E, \Xi_{\lambda})} \big[ \varphi^{-1}(F_j(\varphi)) \big]$$
  where \ $F_j(\varphi)$ \ is the normal sub-module of rank \ $j$ \ of the \ $[\lambda]-$primitive theme \ $\varphi(E)$, with the convention that \ $F_j(\varphi) = \varphi(E)$ \ when the rank of\ $\varphi$ \ is \ $\leq j$.
  \item The ss-depth of \ $E$ \ is equal to \ $d$ \ if and only if \ $d$ \ is the maximal rank of a \ $[\lambda]-$primitive quotient theme of \ $E$.
  \item The ss-depth of \ $E$ \ is equal to \ $d$ \ if and only if \ $d$ \ is the maximal rank of a normal \ $[\lambda]-$primitive sub-theme of \ $E$.
  \end{enumerate} 
   \end{prop}

   \parag{Proof of proposition \ref{ss sequence}} Let us prove i) by induction on \ $j$. As the case \ $j = 1$ \ is obvious, let us assume that \ $j \geq 2$ \ and that the result is proved for \ $j-1$. Let \ $T$ \ a \ $[\lambda]-$primitive theme in \ $E$, and let \ $F_{j-1}(T)$ \ be its normal submodule of rank \ $j-1$ \ (equal to \ $T$ \ if the rank of \ $T$ \ is less than \ $j-1$). Then by the induction hypothesis, we have \ $F_{j-1}(T) \subset S_{j-1}(E)$. Then we have a \ $\A-$linear map\ $T\big/F_{j-1}(T) \rightarrow E\big/S_{j-1}(E)$. If  the rank of \ $T$ \ is at most \ $j$, then \ $T\big/F_{j-1}(T) $ \ has rank at most 1 and its image is in \ $S_1(E\big/S_{j-1}(E))$. So \ $T \subset S_j(E)$.\\
   To prove ii) we also make an induction on \ $j$. The case \ $j = 1$ \ is obvious. So we may assume \ $j \geq 2$ \ and the result proved for \ $j-1$. Let \ $\varphi : S_j(E )\to T$ \ a surjective map on a \ $[\lambda]-$primitive theme \ $T$. By the inductive hypothesis we have \ $\varphi(S_{j-1}(E)) \subset  F_{j-1}(T)$. So we have an induced surjective map
   $$ \tilde{\varphi} : S_j(E)\big/S_{j-1}(E) \to T\big/F_{j-1}(T).$$
  As \ $S_j(E)\big/S_{j-1}(E)$ \ is semi-simple, the image of \ $\tilde{\varphi}$ \ has rank \ $\leq 1$. It shows that \ $T$ \ has rank \ $\leq j$.\\
  To prove iii) consider first a \ $\A-$linear map\ $\varphi : E \to \Xi_{\lambda}$. As \ $\varphi(E)$ \ is a \ $[\lambda]-$primitive theme, \ $\varphi(S_j(E))$ \ is a \ $[\lambda]-$primitive theme quotient of \ $S_j(E)$. So its rank is \ $\leq j$ \ and we have \ $\varphi(S_j(E)) \subset F_j(\varphi)$.\\
  Conversely, for any \ $\A-$linear map\ $\varphi : E \to \Xi_{\lambda}$, the image \ $\varphi(S_j(E))$ \ is a  \ $[\lambda]-$primitive quotient theme of \ $S_j(E)$. So its rank is \ $\leq j$ \ and it is contained in \ $F_j(\varphi)$.\\
  Let us prove iv). If \ $S_d(E) = E$ \ then any \ $[\lambda]-$primitive sub-theme in \ $E$ \ has rank \ $\leq d$ \ thanks to ii). Conversely, assume that for any \ $[\lambda]$ \  any \ $[\lambda]-$primitive sub-theme of \ $E$ \ has rank \ $\leq d-1$ \ and \ $S_{d-1}(E) \not= E$. Then choose a \ $\A-$linear map \ $\varphi : E \to \Xi_{\lambda}$ \ such that \ $\varphi^{-1}(F_{d-1}(\varphi)) \not= E$. Then \ $\varphi(E)$ \ is a \ $[\lambda]-$primitive theme of rank \ $d$ \ which is a quotient of \ $E$, thanks to the following lemma \ref{clef 0}.\\
  The point v) is an easy consequence of the fact that \ $d(E) = d(E \otimes E_{\lambda})$ \ for any \ $\lambda \in \C$ \ and any regular (a,b)-module \ $E$ \ and the fact that for a fresco \ $E$ \ and an integer \ $N$ \ large enough, \ $E^{*}\otimes E_{N}$ \ is again a fresco. So \ $d(E) = d(E^{*}\otimes E_{N})$ \ is the maximal rank of a \ $[N-\lambda]-$primitive quotient theme  \ $T$ \ of \ $E^{*}\otimes E_{N}$ ; for such \ $T$, $T^{*}\otimes E_{N}$ \ is a \ $ [\lambda]-$primitive sub-theme of \ $E$ \ of rank \ $d(E)$. The point  i) allows to conclude.$\hfill \blacksquare$\\
  
  \begin{lemma}\label{clef 0}
Let \ $E$ \ be a rank \ $k$ \ \ $[\lambda]-$primitive theme and denote by \ $F_j$ \ its normal rank \ $j$ \ submodule. Let \ $x \in E \setminus F_{k-1}$. Then the (a,b)-module \ $\A.x \subset E$ \ is a rank \ $k$ \ theme.
\end{lemma}

\parag{Proof}We  may assume \ $E \subset \Xi_{\lambda}^{(k-1)}$ \ and then (see [B.10]) we have the equality  \ $F_{k-1} = E \cap \Xi_{\lambda}^{(k-2)}$. So \ $x$ \ contains a non zero term with \ $(Log\, s)^{k-1}$ \ and then the result is clear. $\hfill \blacksquare$\\

We conclude this paragraph showing that any geometric (a,b)-module contains a fresco with finite \ $\C-$codimension.
  
 \begin{lemma}\label{fresco in geom.}
 Let \ $E$ \ be a geometric (a,b)-module of rank \ $k$. There exists a fresco \ $E' \subset E$ \ with rank \ $k$ ; so the quotient \ $E\big/E'$ \ is a finite dimensional complex vector space.
 \end{lemma}
 
 \parag{proof} We shall prove this fact by induction on the rank \ $k$ \ of  \ $E$.  As the statement is obvious for \ $k \leq 1$, assume \ $k \geq 2$ \ and the result proved in rank \ $k-1$. As there exists a normal rank \ $1$ \ submodule of \ $E$, consider an exact sequence
 $$ 0 \to E_{\lambda} \to E \overset{\pi}{\longrightarrow}  F \to 0 $$
 where \ $F$ \ is a rank \ $k-1$ \ geometric (a,b)-module. Let \ $x \in E$ \ such that \ $\pi(x)$ \ generates a rank \ $k-1$ \ fresco in \ $F$. Let \ $P \in \A$ \ be a monic degree \ $k-1$ \ polynomial in \ $a$ \ with coefficients in \ $\C[[b]]$, which generates the annihilator of \ $\pi(x)$ \ in \ $F$. Then \ $P.x$ \ is in \ $E_{\lambda}$. We may assume that \ there exists an invertible element \ $S \in \C[[b]]$ \ and an integer \ $m$ \ such that \ $S.P.x = b^{m}.e_{\lambda}$, where \ $e_{\lambda}$ \ is a standard generator of \ $E_{\lambda}$, because, if we have \ $P.x = 0$ \ we may replace \ $x$ \ by \ $ x + b^{N}.e_{\lambda}$, and for \ $ N \gg 1$ \ we have \ $P.b^{N}.e_{\lambda} \not= 0$. Then we have \ $\A.x \cap E_{\lambda} = b^{m.}E_{\lambda}$ \ and so the exact sequence
 $$ 0 \to E_{\lambda}\big/b^{m}.E_{\lambda} \to E\big/\A.x \to F\big/\A.\pi(x) \to 0 $$
 gives the finiteness of the complex vector space \ $E\big/\A.x$. $\hfill\blacksquare$\\

\subsection{A characterization of semi-simple frescos.}

We begin by a simple remark : A  \ $[\lambda]-$primitive theme is semi-simple if and only if it has rank \ $\leq 1$. This is an easy consequence of the fact that the saturation of a rank \ $2$ \ $[\lambda]-$primitive theme is one of the rank \ $2$ \ (a,b)-modules considered in the remark following the definition \ref{ss 0} which are not semi-simple

\begin{lemma}\label{semi-simple}
A  geometric (a,b)-module \ $E$ \ is {\bf semi-simple} if and only if any  quotient of \ $E$ \ which is a \ $[\lambda]-$primitive theme for some \ $[\lambda] \in \mathbb{Q}\big/\mathbb{Z}$ \   is of rank \ $\leq 1$.
\end{lemma}

\parag{proof} The condition is clearly necessary as a quotient of a semi-simple (a,b)-module is semi-simple (see corollary \ref{quot. ss 2}), thanks to the remark above. Using the application of the lemma 
\ref{primitive ss}, it is enough to consider the case of a \ $[\lambda]-$primitive\ $E$ \ to prove that the condition is sufficient. Let \ $\varphi : E \to \Xi_{\lambda}^{(N)} \otimes V$ \ be an embedding of \ $E$ \ which exists thanks to the embedding theorem 4.2.1. of [B.09], we obtain that each component of this map in a basis \ $v_1, \dots v_p$ \ of \ $V$ \ has rank at most \ $1$ \ as its image is a \ $[\lambda]-$primitive theme. Then each of these images is isomorphic to some \ $E_{\lambda+q}$ \ for some integer \ $q$. So we have in fact an embedding of \ $E$ \  in a direct sum of \ $E_{\lambda+ q_i}$ \ and \ $E$ \ is semi-simple. $\hfill \blacksquare$

\parag{Remark} The preceeding proof shows that a geometric (a,b)-module  is a non zero  \ $[\lambda]-$primitive theme for some \ $[\lambda] \in \mathbb{Q}\big/\mathbb{Z}$ \ if and only if \ $S_1(E)$ \ has rank  $1$ (compare with the theorem 3.1.7 of [B.10]). $\hfill \square$

\begin{lemma}\label{all J-H}
Let \ $E$ \ be a semi-simple fresco  with rank \ $k$ \ and let \ $\lambda_1, \dots, \lambda_k$ \ be the numbers associated to a J-H. sequence of \ $E$. Let \ $\mu_1, \dots, \mu_k$ \ be a twisted permutation\footnote{This means that the sequence \ $\mu_j+j, j \in [1,k]$ \ is a permutation (in the usual sens) of \ $\lambda_j+j , j\in [1,k]$.} of \ $\lambda_1, \dots, \lambda_k$. Then there exists a J-H. sequence for \ $E$ \ with quotients corresponding to \ $\mu_1, \dots, \mu_k$.
\end{lemma}

\parag{Proof} As the symetric group \ $\mathfrak{S}_k$ \  is generated by the transpositions \ $t_{j,j+1}$ \ for \ $j \in [1,k-1]$, it is enough to show that, if \ $E$ \ has a J-H. sequence with quotients given by the numbers \ $\lambda_1, \dots, \lambda_k$, then there exists a J-H. sequence for \ $E$ \ with quotients \ $\lambda_1, \dots, \lambda_{j-1}, \lambda_{j+1}+1,\lambda_j-1, \lambda_{j+2}, \dots, \lambda_k$ \ for \ $j \in [1,k-1]$. But \ $G : = F_{j+1}\big/F_{j-1}$ \ is a rank 2  sub-quotient of \ $E$ \ with an exact sequence
$$ 0 \to E_{\lambda_j} \to G \to E_{\lambda_{j+1}} \to 0 .$$
As \ $G$ \ is  a rank 2 semi-simple fresco, it admits also an exact sequence
$$ 0 \to G_1 \to G \to G\big/G_1 \to 0$$
with \ $G_1 \simeq E_{\lambda_{j+1}+1}$ \ and \ $G\big/G_1 \simeq E_{\lambda_j-1}$. Let \ $q : F_{j+1} \to G$ \ be the quotient map. Now the J-H. sequence for \ $E$ \ given by
$$ F_1, \dots, F_{j-1}, q^{-1}(G_1), F_{j+1}, \dots ,F_k = E$$
 satisfies our requirement. $\hfill \blacksquare$
 
 \parag{Remark} If \ $E$ \ is a semi-simple geometric (a,b)-module, we may have \\
  $ G \simeq E_{\lambda_{j}}\oplus E_{\lambda_{j+1}}$ \ in the proof above, and then the conclusion does not hold. $\hfill \square$
 
  \begin{prop}\label{crit. ss}
 Let \ $E$ \ be a  \ $[\lambda]-$primitive fresco. A necessary and sufficient condition in order that \ $E$ \ is semi-simple is that it admits a J-H. sequence with quotient corresponding to \ $\mu_1, \dots, \mu_k$ \ such that the sequence \ $\mu_j+j$ \ is strictly decreasing.
 \end{prop}
 
 \parag{Remarks}
 \begin{enumerate}
 \item As a fresco is semi-simple if and only if for each \ $[\lambda]$ \ its  \ $[\lambda]-$primitive  part is semi-simple, this proposition gives also a criterium to semi-simplicity for any fresco.
  \item This criterium is a very efficient tool to produce easily examples of semi-simple frescos.
  \end{enumerate}
 
 \parag{Proof} Remark first that if we have, for a   fresco \ $E$, a J-H. sequence \ $F_j, j \in [1,k]$ \  such that \ $\lambda_j+j = \lambda_{j+1}+j+1$ \ for some \ $j \in [1,k-1]$, then \ $F_{j+1}\big/F_{j-1}$ \ is a sub-quotient of \ $E$ \ which is a  \ $[\lambda]-$primitive theme of rank  2. So \ $E$ \ is not semi-simple. As a consequence, when a  fresco \ $E$ \ is semi-simple the principal J-H. sequence corresponds to a strictly increasing sequence \ $\lambda_j+j$. Now, thanks to the previous lemma we may find a J-H. sequence for \ $E$ \  corresponding to the strictly decreasing order for the sequence \ $\lambda_j+j$. \\

 No let us prove the converse. We shall use the following lemma.
 
 \begin{lemma}\label{utile}
 Let \ $F$ \ be a rank \ $k$ \ semi-simple  \ $[\lambda]-$primitive fresco and let \ $\lambda_j+j$ \ the strictly increasing sequence corresponding to its principal J-H. sequence. Let \ $\mu \in [\lambda]$ \ such that \ $0 < \mu+k < \lambda_1+1$. Then any fresco \ $E$ \ in an exact sequence
 $$ 0 \to F \to E \to E_{\mu} \to 0 $$
 is semi-simple (and \ $[\lambda]-$primitive).
 \end{lemma}
 
 \parag{Proof} The case \ $k=1$ \ is clear from the classification of rank \ $2$ \ frescos; so assume that \ $k \geq 2$ \ and that we have a rank 2 quotient \ $\varphi : E \to T$ \ where \ $T$ \ is a \ $[\lambda]-$primitive  theme. Then \ $Ker\, \varphi\, \cap F$ \ is a normal sub-module of \ $F$ \ of rank \ $k-2$ \ or \ $k-3$ \ (for \ $k \geq 3$) in \ $E$. If \ $Ker\,\varphi \,\cap F$ \ is of rank \ $k-3$, the rank of \ $F\big/(Ker\, \varphi\, \cap F)$ \ is \ $2$ \  and it injects in \ $T$ \ via \ $\varphi$. So \ $F\big/(Ker\, \varphi \cap F)$ \ is a rank 2 \ $[\lambda]-$primitive theme. As it is semi-simple, because \ $F$ \ is semi-simple, we get a contradiction.\\
 So the rank of \ $F\big/(Ker\, \varphi\, \cap F)$ \ is \ $1$ \ and we have an exact sequence
 $$ 0 \to F\big/(Ker\, \varphi\, \cap F) \to T \to E\big/F \to 0 .$$
 Put \ $F\big/(Ker\, \varphi\, \cap F) \simeq F_1(T)  \simeq E_{\nu}$. Because \ $T$ \ is a \ $[\lambda]-$primitive  theme, we have the inequality \ $\nu+1 \leq \mu+2$. Looking at a J-H. sequence of \ $E$ \ ending by
 $$ \dots \subset Ker\, \varphi \cap F \subset  \varphi^{-1}(F_1(T)) = F \subset E $$
 we see that \ $\nu+k-1$ \ is in the set \ $\{ \lambda_j+j, j \in [1,k]\}$ \  and, as \ $\lambda_1+1$ \ is the infimum of this set, we obtain  \ $\lambda_1+1 \leq \nu+k- 1 \leq \mu+k $ \ contradicting our assumption. $\hfill \blacksquare$
 
 \parag{End of proof of the proposition \ref{crit. ss}} Now we shall prove by induction on the rank of a \ $[\lambda]-$primitive fresco \ $E$ \ that if it admits a J-H. sequence corresponding to a strictly decreasing sequence \ $\mu_j +j$, it is semi-simple. As the result is obvious in rank 1, we may assume \ $k \geq 1$ \ and the result proved for \ $k$. So let \ $E$ \ be a fresco of rank \ $k+1$ \ and let \ $F_j, j\in [1,k+1]$ \ a J-H. sequence for \ $E$ \ corresponding to the strictly decreasing sequence \ $\mu_j+j, j \in [1,k+1]$. Put \ $F_{j}\big/F_{j-1} \simeq E_{\mu_j}$ \ for all \ $j \in [1,k+1]$, define \ $F : = F_k$ \ and \ $\mu : = \mu_{k+1}$; then the induction hypothesis gives that \ $F$ \ is semi-simple and we apply the previous lemma to conclude. $\hfill \blacksquare$\\
 
 The following interesting corollary is an obvious consequence of the previous proposition.
  
 \begin{cor}\label{ss base 3}
  Let \ $E$ \ be a   fresco and let \ $\lambda_1, \dots, \lambda_k$ \ be the numbers associated to any J-H. sequence of \ $E$. Let \ $\mu_1, \dots, \mu_d$ \ be the numbers associated to any  J-H. sequence of \ $S_1(E)$. Then, for \ $j \in [1,k]$,  there exists a rank 1 normal sub-module of \ $E$ \ isomorphic to \ $E_{\lambda_j+j-1}$ \ if and only if there exists \ $i \in [1,d]$ \ such that we have  \ $\lambda_j+j-1 = \mu_i+i-1$.
  \end{cor}

 Of course, this gives the list of all isomorphy classes of  rank 1 normal sub-modules  contained in \ $E$. So, using shifted duality, we get also the list of all isomorphy classes of  rank 1 quotients of \ $E$. It is interesting to note that this gives the list of the possible initial exponents for maximal logarithmic terms which appear in the asymtotic expansion of a given relative de Rham cohomology class after integration on any vanishing cycle in the spectral subspace of the monodromy associated to the eigenvalue \ $exp(2i\pi.\lambda)$.\\
 
 \subsection{An interesting example.}
 
 Our aim is to produce a rank \ $4$ \ $[\lambda]-$primitive fresco \ $E$ \ with the following properties:
\begin{enumerate}[i)]
\item It semi-simple part \ $S_{1}(E)$ \  has rank \ $2$ \ and equal to \ $F_{2}(E)$.
\item It co-semi-simple part \ $\Sigma_{1}(E)$ \ has rank \ $1$ \ and equal to \ $F_{1}(E)$.
\item There exists two normal rank \ $2$ \ sub-themes in \ $E$ \ with different fundamental invariants (but, of course, with \ $F_{1}(T) = \Sigma_{1}(E)$).
\end{enumerate}
Of course the dual of \ $E$ \ has a rank \ $3$ \ semi-simple part and a rank \ $2$ \ co-semi-simple part. So the rank of \ $S_{1}(E^{*})$ \ is different from the rank of \ $S_{1}(E)$.\\

This example shows that the semi-simple and co-semi-simple filtrations of a geometric (a,b)-module do not behave as in the case of the semi-simple and co-semi-simple filtrations of a complex vector space associated to an endomorphism.\\

Fix \ $ \lambda_{1} > 3$ \ a rational number and \ $p_{i}, i = 1,2,3$ \ natural integers at least equal to \ $2$ \ and define
$$ \lambda_{2} : = \lambda_{1} + p_{1} -1 \quad  \lambda_{3} : = \lambda_{2} + p_{2} -1 \quad  \lambda_{4} : = \lambda_{3} + p_{3} -1. $$
Fix also two non zero  complex numbers \ $\alpha, \beta$ \ and define
$$ E : = \A\big/\A.P \quad {\rm with} \quad P : = (a - \lambda_{1}.b).S^{-1}.(a - \lambda_{2}.b).(a - \lambda_{3}.b).(a - \lambda_{4}.b) $$
where \ $S : = 1 + \alpha.b^{p_{1}+p_{2}} + \beta.b^{p_{1}+p_{2}+p_{3}}$.

\parag{proof} First remark that \ $E$ \ is not semi-simple because the identity
$$ (a - \lambda_{2}.b).(a - \lambda_{3}.b) = (a - (\lambda_{3}+1).b).(a - (\lambda_{2}-1).b) $$
and the fact that \ $\alpha \not= 0$, show that we have a rank \ $2$ \ (normal) sub-theme in \ $E$ \ with fundamental invariants \ $(\lambda_{1}, \lambda_{3}+1)$ \ and parameter \ $\alpha \not= 0$.\\
As the quotient \ $E\big/F_{1}(E)$ \ is semi-simple, we conclude that \ $\Sigma_{1}(E) = F_{1}(E)$.\\
To prove that \ $S_{1}(E) = F_{2}(E)$, we first remark that the inclusion \ $ F_{2}(E) \subset S_{1}(E)$ \ is clear as \ $F_{2}$ \ is semi-simple. To prove the opposite inclusion we shall show that any rank \ $1$ \ sub-module in \ $E$ \ is contained in \ $F_{2}(E)$. This is enough to conclude because if \ $S_{1}(E)$ \  has rank \ $3$, it has a rank \ $1$ \ sub-module\footnote{for instance the maximal simple pole sub-module \ $S_{1}(E)^{b}$ \ of \ $S_{1}(E)$ \ has a decomposition \ $(F_{2}\cap S_{1}(E)^{b})\oplus L $ \ so \ $L \cap F_{2} = 0$.} \ $L$ \ with \ $L \cap F_{2} = 0$.\\
For that purpose consider the standard basis \ $e_{1}, e_{2}, e_{3}, e_{4}$ \ of \ $E$ \ associated to its principal J.H. sequence. It satisfies
$$ (a - \lambda_{4}.b).e_{4} = e_{3}, \quad (a - \lambda_{3}.b).e_{3} = e_{2}, \quad (a - \lambda_{2}.b).e_{2} = S.e_{1},  \quad (a - \lambda_{1}.b).e_{1} = 0 .$$
Consider now an element in \ $E$, write \ $ x : = U_{4}.e_{4} + U_{3}.e_{3} + U_{2}.e_{2} + U_{1}.e_{1}$ \ and assume that it satisfies : $(a - \mu.b).x = 0 $ \ for some complex number \ $\mu$. This equality is equivalent to the system of equations
\begin{align*}
& b^{2}.U'_{4} + (\lambda_{4} - \mu).b.U_{4} = 0 \\
& U_{4} + b^{2}.U'_{3} + (\lambda_{3} - \mu).b.U_{3} = 0 \\
&  U_{3} + b^{2}.U'_{2} + (\lambda_{2} - \mu).b.U_{2} = 0 \\
& S.U_{2} + b^{2}.U'_{1} + (\lambda_{1} - \mu).b.U_{1} = 0 \\
\end{align*}
If we have \ $ \mu - \lambda_{4} \not\in \mathbb{N}$, the first equation gives \ $U_{4} = 0 $. If \ $ \mu = \lambda_{4} + q$, for some natural integer \ $q$, we find \ $U_{4} = \rho.b^{q}$ \ for some complex number \ $\rho$. In the case \ $U \not= 0$ \ the second equation forces \ $ q \geq 1$ \ and becomes \ $ \rho.b^{q-1} + b.U'_{3} -(p_{3}+q-1).U_{3} = 0 $ \ and so we have
$$U_{3} = \sigma.b^{q-1} + \tau.b^{p_{3}+ q-1} \quad {\rm with} \quad \sigma = \rho\big/p_{3} .$$
Now the third equation forces \ $q \geq 2$ \ and becomes
$$ (\sigma.b^{q-2} + \tau.b^{p_{3}+q -2}) + b.U'_{2} -(p_{2}+p_{3}+q -2).U_{2} = 0 .$$
Then \ $U_{2}$ \ is given by
$$ U_{2} = \sigma'.b^{q-2} + \tau'.b^{p_{3}+q -2} + \eta.b^{p_{2}+p_{3}+q-2} $$
where the complex number \ $\eta$ \ is arbitrary and we have
$$ \sigma' = \sigma\big/(p_{2} + p_{3})          \quad  \tau' =   \tau\big/p_{2} .              $$
Then the last equation will have a solution if and only if  \ $U_{2}(0) = 0$, so \ $q \geq 3$, and there is no term in \ $b^{p_{1}+p_{2}+p_{3}+ q-2}$ \ in \ $S.U_{2}$. As we assume that \ $\alpha.\beta \not= 0$ \ this is the case only when \ $ \sigma' = \tau' = 0$, and so when \ $\sigma = \tau = 0$. Then we have \ $U_{4} = U_{3} = 0$ \ and \ $x$ \ is in \ $F_{2}$.\\
We come back to the case  \ $U_{4} = 0$. Then if \ $\mu - \lambda_{3} \not\in \mathbb{N}$ \ the second equation implies \ $U_{3} = 0$, which conclude this case.\\
So we may assume \ $ \mu = \lambda_{3} + q$. Then we have \ $ U_{3} = \rho.b^{q}$. Again the third equation forces \ $q \geq 1$ \ and becomes
$$ \rho.b^{q-1} + b.U'_{2} -(p_{2} + q -1).U_{2} = 0 .$$
So we have \ $ U_{2} = \rho'.b^{q-1} + \sigma.b^{p_{2} +q -1}$ \ with \ $\rho' = \rho\big/p_{2}$. Then, to solve the last equation, it is necessary that \ $S.U_{2}$ \ has no  constant term  (so \ $q \geq2$) and no term in \ $b^{p_{1}+ p_{2} +q-1}$ ; this implies \ $ \rho' =  \rho = 0$ \ as we assume \ $\alpha \not= 0$. So we must have  \ $U_{3} = 0$ \ and  \ $x$ \ is in \ $F_{2}$ \ in all cases. $\hfill \blacksquare$

\section{ Numerical criteria for semi-simplicity.}

\subsection{Polynomial dependance.}

All \ $\C-$algebras have a unit.\\

When we consider a sequence of algebraically independent variables \ $\rho : = (\rho_i)_{i \in \mathbb{N}}$ \ we shall denote \ $\C[\rho]$ \ the \ $\C-$algebra generated by these variables. Then \ $\C[\rho][[b]]$ \ will be the commutative \ $\C-$algebra of formal power series
$$ \sum_{\nu = 0}^{\infty} \ P_{\nu}(\rho).b^{\nu} $$
where \ $P_{\nu}(\rho)$ \ is an element in \ $\C[\rho]$ \ so a polynomial in \ $\rho_0, \dots, \rho_{N(\nu)} $ \ where \ $N(\nu)$ \ is an integer depending on \ $\nu$. So each coefficient in the formal power serie in \ $b$ \ depends only on a finite number of the variables \ $\rho_i$.

\begin{defn}\label{polyn. depend.}
Let \ $E$ \ be a (a,b)-module and let \ $e(\rho)$ \ be a family of elements in  \ $E$ \ depending on a family of variables \ $(\rho_i)_{i\in \mathbb{N}}$. We say that \ $e(\rho)$ \ {\bf\em depends polynomially on \ $\rho$} \ if there exists a fixed \ $\C[[b]]-$basis \ $e_1, \dots, e_k$ \ of \ $E$ \ such that
$$ e(\rho) = \sum_{j=1}^k \ S_j(\rho).e_j $$
where \ $S_j$ \ is for each \ $j \in [1,k]$ \ an element in the algebra  \ $\C[\rho][[b]]$.
\end{defn}

\parag{Remarks} \begin{enumerate}
\item It is important to note that when \ $e(\rho)$ \ depends polynomially of \ $\rho$, then \ $a.e(\rho)$ \ also. Then for any \ $u \in \A$, again \ $u.e(\rho)$ \ depends polynomially on \ $\rho$. 
\item It is easy to see that we obtain an equivalent condition on the family \ $e(\rho)$ \ by asking that the coefficient of \ $e(\rho)$ \ are in \ $\C[\rho][[b]]$ \ in a \ $\C[[b]]-$basis  of \ $E$ \  whose elements depend polynomially of \ $\rho$.
\item The invertible elements in the algebra \ $\C[\rho][[b]]$ \ are exactly those elements with a constant term in \ $b$ \ invertible in the algebra \ $\C[\rho]$.  As we assume that the variables \ $(\rho)_{i\in \mathbb{N}}$ \ are algebraically independent, the invertible elements are  those  with a constant term in \ $b$ \ in \ $\C^*$.
 $\hfill \square$\\
\end{enumerate}

\begin{prop}\label{gen. dep.}
Let \ $E$ \ be a rank \ $k$ \  fresco and let \ $\lambda_1, \dots, \lambda_k$ \ its fundamental invariants. Let \ $e(\rho)$ \ be a family of generators of \ $E$ \ depending polynomially on a family of algebraically independent variables \ $(\rho_i)_{i \in \mathbb{N}}$. Then there exists \ $S_1, \dots, S_k$ \ in \ $\C[\rho][[b]]$ \ such that \ $S_j(\rho)[0] \equiv 1$ \ for each \ $j \in [1,k]$ \ and such that the annihilator of \ $e(\rho)$ \ in \ $E$ \ is generated by the element of \ $\A$
$$ P(\rho) : = (a - \lambda_1.b).S_1(\rho)^{-1} \dots S_{k-1}(\rho)^{-1}.(a - \lambda_k.b).S_k(\rho)^{-1} .$$
\end{prop}

\parag{Proof} The key result to prove this proposition is the rank \ $1$ \ case. In this case we may consider a standard generator \ $e_1$ \ of \ $E$ \ which is a \ $\C[[b]]-$basis of \ $E$ \ and satisfies
$$ (a - \lambda_1.b).e_1 = 0 .$$
Then, by definition, we may write
$$ e(\rho) = S_1(\rho).e_1 $$
where \ $S_1 $ \ is in \ $\C[\rho][[b]]$ \ is invertible in this algebra, so has a  constant term in \ $\C^*$. Up to normalizing \ $e_1$, we may assume that \ $S_1(\rho)[0] \equiv 1 $ \ and then define 
 $$P(\rho) : = (a - \lambda_1.b).S_1(\rho)^{-1} .$$
It clearly generates the annihilator of \ $e(\rho)$ \ for each \ $\rho$.\\
Assume now that the result is already proved for the rank \ $k-1 \geq 1$. Then consider the family \ $[e(\rho)]$ \ in the quotient \ $E\big/F_{k-1}$ \ where \ $F_{k-1}$ \ is the rank \ $k-1$ \ sub-module of \ $E$ \ in its principal J-H. sequence. Remark first that \ $[e(\rho)]$ \ is a family of generators of
 \ $E\big/F_{k-1}$ \ which depends polynomially on \ $\rho$. This is a trivial consequence of the fact that we may choose a \ $\C[[b]]-$basis \ $e_1, \dots, e_k$ \  in \ $E$ \ such that \ $e_1, \dots, e_{k-1}$ \ is a \ $\C[[b]]$ \ basis of \ $F_{k-1}$ \ and \ $e_k$ \ maps to a standard generator of \ $E\big/F_{k-1}\simeq E_{\lambda_k}$. Then the rank \ $1$ \  case gives \ $S_k \in \C[\rho][[b]]$ \ with \ $S_k(\rho)[0] = 1$ \  and such that \ $(a - \lambda_k.b).S_k(\rho)^{-1}.e(\rho)$ \ is in \ $F_{k-1}$ \ for each \ $\rho$. But then, thanks to the corollary \ref{gen. Fj}, it is a family of generators of \ $F_{k-1}$ \ which depends polynomially on \ $\rho$ \ and the inductive assumption allows to conclude. $\hfill \blacksquare$\\

Fix now the fundamental invariants \ $\lambda_1, \dots, \lambda_k$ \ of some fresco. 

\begin{defn}\label{polyn. depend. E}
Consider a complex valued  function \ $f$ \  defined on a subset \ $\mathcal{F}_0$ \ of the isomorphism classes \ $\mathcal{F}(\lambda_1, \dots, \lambda_k)$ \ of  frescos with fundamental invariants \ $\lambda_1, \dots, \lambda_k$. We shall say that  \ $f$ \  {\bf depends polynomially on the isomorphism class  \ $[E] \in \mathcal{F}_0 $ \ of the fresco \ $E$} \ (or simply that \ $f$ \ is a polynomial on \ $\mathcal{F}_0$) if the following condition is satisfied :\\
Let \ $s$ \ be the collection of algebraically independent variables corresponding to the non constant  coefficients of \ $k$ \ \'elements \ $S_1, \dots, S_k$ \ in \ $\C[[b]]$ \ such that \\
 $S_j(0) = 1\quad \forall j \in [1,k]$ \ and consider for each value of \ $s$ \ the rank \ $k$ \   fresco \ $E(s) : = \A \big/ \A.P(s) $ \ where 
 $$P(s) : =  (a - \lambda_1.b)S_1(s)^{-1} \dots (a - \lambda_k.b).S_k(s)^{-1} $$
 where \ $S_1(s), \dots, S_k(s)$ \ correspond to the given values for \ $s$.\\
 Then there exists a polynomial \ $F \in \C[s]$ \ such that for each value of  \ $s$ \ such that \ $[E(s)]$ \ is in \ $\mathcal{F}_0$, the value of \ $F(s)$ \ is equal to\ $f([E(s)])$.
 \end{defn}
 
 \begin{defn}\label{polyn. str,}
 We shall say that a subset in \ $\mathcal{F}(\lambda_1, \dots, \lambda_k)$ \ is {\bf algebraic} if it is the common zero set of a  finite set of polynomial functions on \ $\mathcal{F}(\lambda_1, \dots, \lambda_k)$.
 An algebraic stratification  of \ $\mathcal{F}(\lambda_1, \dots, \lambda_k)$ \ is a finite family  of algebraic subsets
 $$ \emptyset \subset S_N \subset \dots \subset S_0 = \mathcal{F}(\lambda_1, \dots, \lambda_k).$$

 \end{defn}
 
 \parag{Remark} Let \ $\mathcal{G}$ \ be an algebraic subset in \ $\mathcal{F}(\lambda_1, \dots, \lambda_k)$ \ and \ $f : \mathcal{G} \to \C$ \ a polynomial, then \ $\{ f = 0\} \subset \mathcal{G}$ \ is  an algebraic subset of \ $\mathcal{F}(\lambda_1, \dots, \lambda_k)$. $\hfill \square$\\

\begin{defn}\label{polyn. appl.}
Let \ $ f : \mathcal{F}(\lambda_1, \dots, \lambda_k)  \longrightarrow \mathcal{F}(\mu_1, \dots, \mu_l)$ \ an application. We shall say that \ $f$ \ is algebraic if for each polynomial  \ $p : \mathcal{F}(\mu_1, \dots, \mu_l) \to \C$ \ the composition  \ $p\circ f$ \ is a polynomial.
\end{defn}

\parag{Remark} In the situation of the previous definition, let \ $s$ \ the collection of algebraically independent variables corresponding to the non constant  coefficients of \ $k$ \ \'elements \ $S_1, \dots, S_k$ \ in \ $\C[[b]]$ \ satisfying \ $S_j(0) = 1\quad \forall j \in [1,k]$; assume that there exist \ $T_1, \dots, T_l \in \C[s][[b]]$ \ such that for each \ $s$ \  we have $$f([E(s)]) = [\A\big/\A.(a - \mu_1.b).T_1(s)^{-1}\dots T_{l-1}(s)^{-1}.(a - \mu_l.b).T_{l}(s)^{-1}]$$
 where we define \ $E(s) : = \A \big/ \A.P(s) $, then the map \ $f$ \ is algebraic.

\parag{Examples}
\begin{enumerate}
\item Let \ $(s_i)_{i \in \mathbb{N}^*}$ \ a family of algebraically independent variables and define   \\ $S(s) : = 1 + \sum_{i=1}^{\infty} \ s_i.b^i \in \C[s][[b]]$. Let  
 $$E(s) : = \A\big/\A.(a - \lambda_1.b).S(s)^{-1}.(a - \lambda_2.b) $$
 where \ $\lambda_1 > 1$ \ is rational and \ $\lambda_2 : = \lambda_1 + p_1 -1$. For \ $p_1 \in \mathbb{N}^{*}$ \ define \ $\alpha(s) : = s_{p_1}$.  For \ $ p_{1} = 0$ \  define \ $\alpha(s) \equiv 1$ \ and for \ $p_{1}\not\in \mathbb{N}$ \ define \ $ \ \alpha(s) \equiv 0$. Then the number \ $\alpha(s)$ \ depends only of the isomorphism class of the fresco \ $E(s)$ \ and defines a polynomial on \ $\mathcal{F}(\lambda_1, \lambda_2)$.
\item Fix \ $\lambda_1, \dots, \lambda_k$ \ and \ $(i,h) $ \ such that \ $[i, i+h] \subset [1, k]$ \ and consider the map
$$\pi : \mathcal{F}(\lambda_1, \dots, \lambda_k) \longrightarrow \mathcal{F}(\lambda_i, \dots, \lambda_{i+h}) $$
given by \ $\pi([E]) = [F_{i+h}\big/F_{i-1}]$ \ where \ $(F_j)_{j \in [0,k]}$ \ is the principal J.H. sequence of \ $E$. Then this map is obviously algebraic, thanks to the previous remark.
\item Combining the two examples above, we see that for each \ $j \in [1, k-1]$ \ the function
$$ \alpha_j : [E] \mapsto \alpha([F_{j+1}\big/F_{j-1}]) $$
is a polynomial on \ $\mathcal{F}(\lambda_1, \dots, \lambda_k)$.
\item We shall define 
 $$\mathcal{F}_2(\lambda_1, \dots, \lambda_k) : = (\alpha_1, \dots, \alpha_{k-1})^{-1}(0) .$$
 This is an algebraic subset of \ $\mathcal{F}(\lambda_1, \dots, \lambda_k)$. Note that if there exists at least  one \ $j \in [1,k-1]$ \ with \ $p_j : = \lambda_{j+1} - \lambda_j + 1 = 0$, then the subset \ $\mathcal{F}_2(\lambda_1, \dots, \lambda_k)$ \ is empty. Note also that in this case there is no semi-simple fresco with fundamental invariants \ $\lambda_1, \dots, \lambda_k$ \ as \ $F_{j+1}\big/F_{j-1}$ \ is a theme.\\
 Conversely, it is easy to see (for instance using the proposition \ref{crit. ss} and the lemma \ref{primitive ss}) that if we have \ $p_j \not= 0$ \ for all \ $ j \in [1,k-1]$ \ there exists a semi-simple rank \ $k$ \  fresco with fundamental invariants \ $\lambda_1, \dots, \lambda_k$.\\
 Remark that a class \ $[E]$ \ is in  the subset \ $\mathcal{F}_2(\lambda_1, \dots, \lambda_k)$ \ if and only if any rank \ $2$ \ sub-quotient of the principal J.H. sequence of \ $E$ \ is semi-simple. $\hfill \square$
 \end{enumerate}

\subsection{The \ $\beta-$invariant.}

In this section we shall fix \ $k \geq 3$ \ and \ $\lambda_1, \dots, \lambda_k$ \ the fundamental invariants of a rank \ $k$ \ $[\lambda]-$primitive semi-simple fresco. So we have \ $p_j \geq 1$ \ for all \ $j \in [1,k-1]$ \ and \ $\lambda_1 > k-1$ \ is a rational number.\\
Let \ $\mathcal{F}_{k-1}(\lambda_1, \dots, \lambda_k)$ \ the subset of \ $\mathcal{F}(\lambda_1, \dots, \lambda_k)$ \ of isomorphism class of rank \ $k$ \ $[\lambda]-$primitive fresco \ $E$ \  with invariants \ $\lambda_1, \dots, \lambda_k$ \ such that \ $F_{k-1}$ \ and \ $E\big/F_1$ \ are semi-simple where \ $(F_j)_{j\in [1,k]} $ \ is  the principal J-H. sequence of \ $E$.\\
Fix \ $E$ \ with \ $[E]$ \ in \ $\mathcal{F}_{k-1}(\lambda_1, \dots, \lambda_k)$. Then \ $E$ \ admits a generator \ $e$ \ such that \ $(a - \lambda_{k-1}.b).(a - \lambda_k.b).e$ \ lies in \ $F_{k-2}$. This is a consequence of the fact that \ $E\big/F_{k-2}$ \ is semi-simple.

\begin{lemma}\label{tech.1}
In the above situation,  let \ $e$ \ be a generator of \ $E$ \ such that 
 $$(a - \lambda_{k-1}.b).(a - \lambda_k.b).e$$
  lies in \ $F_{k-2}$. Then the sub-module \ $G_{k-1}$ \ of \ $E$ \ generated by the element  \\
   $(a - (\lambda_{k-1}-1).b).e$ \ is a normal rank \ $k-1$ \ sub-module of \ $E$ \ which contains \ $F_{k-2}$ \ and which  is in \ $\mathcal{F}_{k-2}(\lambda_1, \dots, \lambda_{k-2}, \lambda_k+1)$.
\end{lemma}

\parag{proof}  Using the identity in \ $\A$ :
 $$(a - \lambda_{k-1}.b).(a - \lambda_k.b) = (a - (\lambda_k+1).b).(a - (\lambda_{k-1}-1).b) $$
 we see that \ $G_{k-1}$ \ contains \ $F_{k-2}$ \ (because \ $(a - \lambda_{k-1}.b).(a - \lambda_k.b).e$ \ is a generator of \ $F_{k-2}$ \ thanks to the corollary \ref{gen. Fj}), that \ $G_{k-1}\big/F_{k-2} \simeq E_{\lambda_k+1}$ \  and so \ $G_{k-1}$ \ admits 
 $$ F_1 \subset \dots \subset F_{k-2} \subset G_{k-1} $$
 as principal J-H. sequence. Then the  fundamental invariants for \ $G_{k-1}$ \ are equal to \ $\lambda_1, \dots, \lambda_{k-2},\lambda_k+1$, and the corank \ $1$ \  term \ $F_{k-2}$ \ is semi-simple. As \ $G_{k-1}\big/F_1$ \ is a sub-module of \ $E\big/F_1$ \ which is semi-simple by assumption, it is semi-simple and we have proved that 
 \ $ [G_{k-1}]$ \ is in \ $ \mathcal{F}_{k-2}(\lambda_1, \dots, \lambda_{k-2}, \lambda_k+1). \hfill \blacksquare$\\

\begin{lemma}\label{tech.2}
 If \ $\varepsilon$ \ is another generator of \ $E$ \ such that \ $(a - \lambda_{k-1}.b).(a - \lambda_k.b).\varepsilon $ \ is in \ $F_{k-2}$, we have
$$ \varepsilon = \rho.e + \sigma.b^{p_{k-1}-1}.(a - \lambda_k.b).e \quad {\rm modulo} \ F_{k-2} $$
where \ $(\rho, \sigma)$ \ is in \ $\C^*\times \C$. And conversely, for any \ $(\rho, \sigma)$ \  in \ $\C^*\times \C$, an element \ $x = \rho.e + \sigma.b^{p_{k-1}-1}.(a - \lambda_k.b).e \quad {\rm modulo} \ F_{k-2} $ \ is a generator of \ $E$ \ and satifies \ $(a - \lambda_{k-1}.b).(a - \lambda_k.b).x$ \ is in \ $F_{k-2}$.
\end{lemma}

\parag{proof} Write \ $\varepsilon = U.e_k + V.e_{k-1} \quad {\rm modulo} \ F_{k-2}$ \ where we defined \ $e_k : = e$ \ and \ $e_{k-1} : = (a - \lambda_k.b).e $ \ and where \ $U, V $ \ are in \ $\C[[b]]$. Now 
$$ (a - \lambda_k.b).\varepsilon = b^2.U'.e_k + U.e_{k-1} + (\lambda_{k-1} - \lambda_k).b.V.e_{k-1} +  b^2.V'.e_{k-1} \quad {\rm modulo} \ F_{k-2} $$
and, as \ $(a - \lambda_{k-1}.b).b^2.U'e_k \in F_{k-1}$ \ implies \ $U' = 0$ \ ( because \ $p_{k-1} + 1 > 0$)   we have  \ $U = \rho \in \C^*$, and we obtain
$$ (a - \lambda_k.b).\varepsilon = \left[ \rho + b^2.V' - (p_{k-1}-1).b.V\right].e_{k-1} \quad {\rm modulo} \ F_{k-2} .$$
If \ $T : = \rho + b^2.V' - (p_{k-1}-1).b.V$ \ we have
$$ (a - \lambda_{k-1}.b).(a - \lambda_k.b).\varepsilon = b^2.T'.e_{k-1}  \quad {\rm modulo} \ F_{k-2} $$
and so  \ $T $ \ is a constant and equal to \ $\rho$. Then \ $V = \sigma.b^{p_{k-1}-1}$ \ for some complex number \ $\sigma$. \\
The converse is easy, thanks to the lemma \ref{gen.} and the previous computation.$\hfill \blacksquare$\\

\begin{cor}\label{tech.3}
Let \ $ \varepsilon = \rho.e + \sigma.b^{p_{k-1}-1}.(a - \lambda_k.b).e \quad {\rm modulo} \ F_{k-2} $ \ be a generator of \ $E$ \ where \ $(\rho, \sigma)$ \ is in \ $\C^*\times \C$. Then the sub-module \ $G : =  \A.(a - (\lambda_{k-1}-1).b).\varepsilon$ \ is equal to the sub-module \ $G^{\tau} : = \A.(a - (\lambda_{k-1}-1).b).\varepsilon(\tau)$ \ where we define \ $ \varepsilon(\tau) : = e_k + \tau.b^{p_{k-1}-1}.e_{k-1}$ \ and \ $\tau : = \sigma/\rho$.
\end{cor}

\parag{proof} Using the lemma \ref{gen.}, we see that \ $\varepsilon$ \ and \ $\varepsilon(\tau)$ \ are generators of \ $E$ \ such that there images by \ $(a - \lambda_{k-1}.b).(a - \lambda_{k}.b)$ \ are in\ $F_{k-2}$. As the sub-module \ $G$ \ contains \ $F_{k-2}$ \ it contains \ $\varepsilon(\tau)$ \ and so \ $G$ \ contains \ $G^{\tau}$. But they are normal and have same rank ; so they are equal.$\hfill \blacksquare$\\

\begin{prop}\label{tech.4}
Assume that we already have proved  that for any fixed fundamental invariants of rank  \ $l \leq (k-1)$ \ frescos the following assertion :
\begin{itemize}
\item There exists a polynomial  \ $\beta : \mathcal{F}_{l-1} \to \C$ \ such that \ $[G] \in \mathcal{F}_{l-1}$ \ is semi-simple if and only if \ $\beta([G]) = 0 $.
\end{itemize}
Then for a given \ $[E] \in \mathcal{F}_{k-1}(\lambda_1, \dots, \lambda_k)$ \  the value of \ $\beta(G^{\tau})$ \ is independent of the choice of \ $\tau$ \ and defines a polynomial on \ $\mathcal{F}_{k-1}(\lambda_1, \dots, \lambda_k)$ \ such that a class \ $[E] \in \mathcal{F}_{k-1}(\lambda_1, \dots, \lambda_k)$ \ is semi-simple if and only if \ $\beta([E]) = 0$.
\end{prop}

\parag{proof} Remark that for \ $k = 2$ \ we have \ $\mathcal{F}_1(\lambda_1, \lambda_2) = \mathcal{F}(\lambda_1, \lambda_2)$ \ and that the polynomial \ $\beta$ \ given by the parameter of the rank \ $2$ \ fresco satisfies our requirement; explicitely this means that :
\begin{enumerate}
\item When \ $\lambda_2 - \lambda_1 + 1 \not\in \mathbb{N} $ \ define \ $\beta \equiv 0 $.
\item When \ $\lambda_2 - \lambda_1 + 1= p \in \mathbb{N} $ \ define \ $\beta([E]) : = \alpha(E) $ \ (recall that \ $\alpha(E)$ \ is the parameter of \ $E$ ; see example 1 at the end of paragraph 3.1.) 
\end{enumerate}
So we may assume now that \ $k \geq 3$. We shall prove that the function \ $\tau \mapsto \beta(G^{\tau})$ \ is a polynomial, using the induction hypothesis. For that purpose, it is enough to show that there exists, for a given \ $[E] \in \mathcal{F}_{k-1}(\lambda_1, \dots, \lambda_k)$, elements \ $T_1, \dots, T_{k-2}$ \ in \ $ \C[\tau][[b]]$ \ such that
$$ T_j(0) = 1 \quad \forall j \in [1,k-2] \quad {\rm and} \quad  G^{\tau} \simeq \tilde{A}\big/\tilde{A}.Q(\tau)$$
where \ 
$$Q_0(\tau) : = (a - \lambda_1.b).T_1(\tau)^{-1}.(a - \lambda_2.b) \dots (a - \lambda_{k-2}.b).T_{k-2}(\tau)^{-1} $$
 and \ $ Q(\tau) : = Q_0(\tau).(a - (\lambda_k+1).b) .$ \\
Now, as \ $\varepsilon(\tau)$ \ depends polynomially on \ $\tau$, the element of \ $F_{k-2}$
$$ (a - \lambda_{k-1}.b).(a - \lambda_k.b).\varepsilon(\tau) $$
which is a generator of \ $F_{k-2}$ \ thanks to the corollary \ref{gen. Fj}, depends polynomially on \ $\tau$. Thanks to the proposition \ref{gen. dep.} we may find \ $T_1, \dots, T_{k-2} \in \C[\tau][[b]]$ \ such that the annihilator of \ $(a - \lambda_{k-1}.b).(a - \lambda_k.b).\varepsilon(\tau) $ \ is equal to \ $\tilde{A}.Q_0(\tau)$. Then we have for each \ $\tau$ \ an isomorphism \ $G^{\tau} \simeq \tilde{A}\big/\tilde{A}.Q(\tau)$, where \ $Q(\tau) : = Q_0(\tau).(a - (\lambda_k+1).b)$.\\
Assume that for some value of \ $\tau$ \ we have \ $\beta(G^{\tau}) = 0$. Then \ $G^{\tau}$ \ is semi-simple, and then \ $F_{k-1} + G^{\tau}$ \ is semi-simple. So its normalization in \ $E$ \ is semi-simple. But this normalization is \ $E$, and so any sub-module of \ $E$ \ is semi-simple. This implies that either we have \ $\beta(G^{\tau}) \not= 0$ \ for each \ $\tau$ \ or \ $\beta(G^{\tau}) = 0$ \ for each \ $\tau$. This means that the polynomial \ $\tau \mapsto \beta(G^{\tau})$ \ is constant, and that this constant is \ $0$ \ if and only if \ $E$ \ is semi-simple. To conclude, we define the function \ $ \beta : \mathcal{F}_{k-1}(\lambda_1, \dots, \lambda_k) \to \C$ \ by \ $\beta(E) : = \beta(G^{0})$. This is a polynomial because the application \ $ [E] \mapsto [G^{0}]$ \ is algebraic with value in \ $\mathcal{F}_{k-2}(\lambda_1, \dots, \lambda_{k-2}, \lambda_{k}+1)$ \ and  thanks to our inductive assumption. $\hfill \blacksquare$\\

Let me give in rank \ $3$ \ a polynomial \ $F \in \C[s]$ \ such that we have \ $F(s) = \beta([E(s)])$ \ for those \ $s$ \ for which \ $E(s)$ \ is in \ $\mathcal{F}_2(\lambda_1,\lambda_2,\lambda_3)$.\\
Assume \ $p_1\geq 1 $ and \ $p_2 \geq 1$. Then the necessary and sufficient condition to be in \ $\mathcal{F}_2(\lambda_1,\lambda_2,\lambda_3)$ \ is \ $s^1_{p_1} = s^2_{p_2} = 0 $.

\begin{lemma}
Assume that  \ $p_1\geq 1 $ and \ $p_2 \geq 1$ \ and  let 
 $$E : = \A\big/\A.(a - \lambda_1.b).S_1^{-1}.(a - \lambda_2.b).S_2^{-1}.(a - \lambda_3.b)$$
with  \ $s^1_{p_1} = s^2_{p_2} = 0 $;  the complex number \ $ \beta(E) $ \ is the coefficient of \ $b^{p_1+p_2}$ \ in \ $V.S_1$ \ where \ $V \in \C[[b]]$ \ is any solution of the differential equation 
 $$ b.V' = p_2.(V - S_2) .$$
\end{lemma}

The proof is left to the reader.\\

This gives the following formula for \ $\beta$ :
$$\beta(s) =  p_2.\sum_{j\not= p_2, j=0}^{p_1+p_2} \ s^1_{p_1+p_2-j}.\frac{s^2_{j}}{p_2-j}.  $$
Remark that for \ $S_2 = 1$ \ we find that \ $F(s)$ \ reduces to \ $s^1_{p_1+p_2}$. Using the commutation relation \ $(a - \lambda_2.b).(a - \lambda_3.b) = (a - (\lambda_3+1).b).(a - (\lambda_2-1).b)$ \ we find in this case that \ $E(s)$ \ has a normal rank \ $2$ \ sub-module with fundamental invariants \ $\lambda_1, \lambda_3+1$ \ and parameter \ $s^1_{p_1+p_2}$ \ which is precisely \ $\alpha(E(s))$.\\

Now we shall show that in the previous situation when the \ $\beta-$invariant of a fresco \ 
$E \in \mathcal{F}_{k-1}(\lambda_1, \dots, \lambda_k)$ \  is not zero it is the parameter of any normal rank \ $2$ \ sub-theme of \ $E$. Although such a rank \ $2$ \ normal sub-theme is not unique (in general), its isomorphism class is uniquely determined from the fundamental invariants of \ $E$ and from the \ $\beta-$invariant of \ $E$. It will be enough to consider the \ $[\lambda]-$primitive case thanks to the following remark.

\parag{Remark} If \ $E \in \mathcal{F}_{k-1}(\lambda_1, \dots, \lambda_k)$ \ satisfies \ $\beta(E) \not= 0$ \ then \ $E$ \ is \ $[\lambda]-$primitive: if \ $\lambda_1 \not= \lambda_k \ modulo \  \mathbb{Z}$, let
$$ i : = \sup\{ j \in [1,k] \ / \  \lambda_j = \lambda_1 \ modulo \ \mathbb{Z} \} < k  .$$
Then the lemma \ref{gen. Fj} applied to the exact sequence of frescos
$$ 0 \to F_i \to E \to E\big/F_i \to 0 $$ 
shows that \ $E$ \ is semi-simple, contradicting our assertion \ $\beta(E) \not= 0$. So \ $i = k$ \ and \ $E$ \ is \ $[\lambda]-$primitive. $\hfill \square$\\

\begin{prop}\label{sub-themes}
Let \ $E$ \ be a rank \ $k \geq 2$ \  fresco such that \ $F_{k-1}$ \ and \ $E\big/F_1$ \ are semi-simple and with \ $\beta(E) \not= 0$. Note \ $\lambda_1, \dots, \lambda_k$ \ the fundamental invariants\footnote{Our hypothesis implies that \ $E$ \  is \ $[\lambda]-$primitive, thanks to the remark above, so \\ $p_{j }: = \lambda_{j+1} -\lambda_{j}+ 1$ \ is a positive integer for each \ $j \in [1,k-1]$.} of \ $E$ \ and \ $p(E) : = \sum_{j=1}^{p-1} \ p_j$. Then there exists at least one rank \ $2$ \ normal  sub-theme in \ $E$ \ and  each rank \ $2$ \ normal sub-theme of \ $E$ \ is  isomorphic to
\begin{equation*}
 \A\big/\A.(a - \lambda_1.b).(1 + \beta(E).b^{p(E)})^{-1}.(a - (\lambda_k+k-2).b). \tag{@}
 \end{equation*}
\end{prop}

\parag{proof} The case \ $k =2$ \ is clear, so we may assume that \ $k \geq 3$ \ and we shall prove the proposition by induction on \ $k$. So assume that the proposition is known for the rank \ $k-1 \geq 2$ \ and let \ $E$ \ be a rank \ $k$ \ ($[\lambda]-$primitive) fresco satisfying our assumptions. The fact that there exists a rank \ $2$ \ normal sub-theme is consequence of the induction hypothesis as  \ $G^{\tau}$ \ for any \ $\tau \in \C$ \ is normal rank \ $k-1$ \ in \ $E$ \ and satisfies again our assumptions.\\
Recall that the fundamental invariants of \ $G^{\tau}$ \ are \ $\mu_1, \dots , \mu_{k-1}$ \ with \ $\mu_j = \lambda_j$ \ for \ $j \in [1,k-2]$ \ and \ $\mu_{k-1} = \lambda_k + 1$. We have also \ $\beta(G^{\tau}) = \beta(E)$ \ for each \ $\tau$, thanks to the proof of the proposition \ref{tech.4}.  As we have  \ $p(G^{\tau}) = p(E)$ \ for each \ $\tau$, and \ $\mu_{k-1} = \mu_1 + p(E) - 1$ \ the inductive hypothesis  implies that any rank \ $2$ \ normal sub-theme of any \ $G^{\tau}$ \ is isomorphic to \ $(@)$. Then, to complete the proof, it is enough to show that any rank \ $2$ \ normal sub-theme of \ $E$ \ is contained in some \ $G^{\tau}$.
\smallskip

Let \ $T$ \ be a rank \ $2$ \ normal sub-theme of \ $E$. We shall first prove that its fundamental invariants are equal to \ $\lambda_1, \lambda_k+k-2$. As \ $E\big/F_1$ \ is semi-simple, the image of \ $T$ \ by the quotient map \ $E \to E\big/F_1$ \ has rank \ $1$ \ and this implies that \ $F_1\cap T$ \ is a rank \ $1$ \ normal sub-module. Then this implies that \ $F_1$ \ is the unique normal rank \ $1$ \ submodule of \ $T$, proving that \ $\lambda_1(T) = \lambda_1$.
\smallskip
We shall prove now that \ $\lambda_2(T) = \lambda_k+k-2$.\\
First remark that the uniqueness of the principal J-H. sequence of \ $E$ \ implies the uniqueness of the quotient map \ $E \to E\big/F_{k-1} \simeq E_{\lambda_k}$ \ because any surjective map \ $q : E \to E_{\lambda_k}$ \ will produce  a principal  a J-H. sequence for \ $E$ \ by adjoining a principal J-H. sequence for \ $Ker\, q$. So a quotient \ $E\big/H$ \ admits a surjective map on \ $E_{\lambda_k}$ \ if  and only if \ $H \subset F_{k-1}$. So \ $H$ \ has to be semi-simple. Then the quotient \ $E\big/T$ \ has no surjective map on \ $E_{\lambda_k}$.\\
The exact sequence of frescos
$$ 0 \to T\big/F_1 \to E\big/F_1 \to E\big/T \to 0  $$
gives the equality  in \ $\A$, thanks to the remark following the proposition \ref{Bernst. fresco} :
 $$P_{T\big/F_1}.P_{E\big/T} = P_{ E\big/F_1} = (a - \lambda_2.b) \dots (a - \lambda_k.b) .$$
 But as \ $E\big/T$ \ is semi-simple and does not have a quotient isomorphic to \ $E_{\lambda_k}$ \ this implies \ $T\big/F_1 \simeq E_{\lambda_k+ k-2}$, proving our claim.\\
 So \ $T$ \ is isomorphic to \ $\A\big/\A.(a - \lambda_1.b).(1 + \beta.b^{p(E)})^{-1}.(a - (\lambda_k+k-2).b)$ \ for some \ $\beta \in \C^*$. Let \ $x \in T$ \ be a generator of \ $T$ \ which is annihilated by 
$$(a - \lambda_1.b).(1 + \beta.b^{p(E)})^{-1}.(a - (\lambda_k+k-2).b).$$
 It satisfies 
\begin{equation*}
(a - (\lambda_k+k-2).b).x \in F_1 \subset F_{k-2} .\tag{*}
\end{equation*}
We shall determine all \ $x \in E$ \ which satisfies the condition \ $(^{*})$ \ modulo \ $F_{k-2}$. Fix a generator \ $e : = e_k$ \ of \ $E$ \ such that \ $e_{k-1} : = (a - \lambda_k.b).e_k$ \ is in \ $F_{k-1}$ \ and satisfies \ $(a - \lambda_{k-1}.b).e_{k-1} \in F_{k-2}$. Then we look for \ $U, V \in \C[[b]]$ \ such that
$$ x : = U.e_k + V.e_{k-1} \quad {\rm satisfies} \quad (a - (\lambda_k+k-2).b).x \in F_{k-2} .$$
This leads to the equations
\begin{align*}
& b^.U' - (k-2).b.U = 0 \\
& U + b^2.V' - (p_{k-1}+k - 3).b.V = 0 
\end{align*}
and so we get
\begin{equation*}
U = \rho.b^{k-2} \quad {\rm and} \quad V = \frac{\rho}{p_{k-1}}.b^{k-3} + \sigma.b^{p_{k-1}+k-3}
 \end{equation*}
 Note that for \ $\rho = 0$ \ we would have \ $x \in F_{k-1}$ \ and this is not possible for the generator of a rank \ $2$ \ theme as \ $F_{k-1}$ \ is semi-simple. Then assuming \ $\rho \not= 0$ \ we may write
 $$ x = \frac{\rho.b^{k-3}}{p_{k-1}}.\left[ p_{k-1}.b.e_k + e_{k-1} + \frac{p_{k-1}.\sigma}{\rho}.b^{p_{k-1}}.e_{k-1}\right] \quad {\rm modulo} \ \  F_{k-2} .$$
 Now recall that the generator of \ $G^{\tau}$ \ is given by
 $$ (a - (\lambda_{k-1}-1).b).\varepsilon(\tau) $$
 where \ $\varepsilon(\tau) = e_k + \tau.b^{p_{k-1}-1}.e_{k-1} $. A simple computation gives
 $$  (a - (\lambda_{k-1}-1).b).\varepsilon(\tau) = p_{k-1}.b.e_k + e_{k-1} + p_{k-1}.\tau.b^{p_{k-1}}.e_{k-1} \quad {\rm modulo} \ F_{k-2} .$$
 As we know that each \ $G^{\tau}$ \ contains \ $F_{k-2}$, we conclude that any such \ $x$ \ is in  \ $G^{\tau}$ \ for \ $\tau = \sigma\big/\rho $, concluding the proof, thanks to the induction hypothesis and the equalities \ $\beta(E) = \beta(G^{\tau})$ \ and \ $p(E) = p(G^{\tau})$.  $\hfill \blacksquare$\\

Using duality we may deduce from this result the following corollary.

\begin{cor}\label{quot. rk 2}
In the situation of the previous proposition \ $E$ \ has a rank \ $2$ \ quotient theme and any such rank \ $2$ \  quotient theme is isomorphic to
$$ \A\big/\A.(a - (\lambda_1-k+2).b).(1 + \beta^*(E).b^{p(E)})^{-1}.(a - \lambda_k.b) $$
where
$$ \beta^*(E) : = (-1)^k\frac{p_1.(p_1+p_2)\dots (p_1+ \dots+p_{k-2})}{p_{k-1}.(p_{k-2}+p_{k-1})\dots (p_2+ \dots p_{k-1})} \beta(E)$$
\end{cor}

The proof is left as an exercice for the reader who may use the following remark.

\parag{Remark} Let \ $E$ \ be a rank \ $2$ \ \ $[\lambda]-$primitive theme with fundamental invariants \ $\lambda_1, \lambda_2$ \ and parameter \ $\alpha(E)$. For \ $\delta \in \mathbb{N}, \delta \gg 1$ \ $E^*\otimes E_{\delta}$ \ is a \ $2$ \ $[\delta-\lambda]-$primitive theme with fundamental invariants\ $(\delta - \lambda_2, \delta - \lambda_1)$ \ and parameter \ $-\alpha(E)$.\\
Then in the situation of the proposition \ref{sub-themes} \ $E^*\otimes E_{\delta}, \delta \gg 1$ \ satisfies again our assumptions and we have \ $\beta(E^*\otimes E_{\delta}) = \beta^*(E)$. $\hfill \square$\\

Using the \ $\beta-$invariant constructed above, we obtain the following result.

\begin{thm}\label{strat. ss}
  Fix a positive  integer \ $k$ \ and  let  \ $\lambda_1, \dots,  \lambda_k$ the fundamental invariants of some rank \ $k$ \ fresco. Then there exists a natural algebraic stratification \ $\mathcal{F}_h, h \in [1,k]$ \ of the subset \ $\mathcal{F}(\lambda_1, \dots, \lambda_k) $ \  with the following properties:
  \begin{enumerate}[i)]
 
 \item A class \ $[E] \in \mathcal{F}(\lambda_1, \dots, \lambda_k) $ \ is in \ $\mathcal{F}_{h}$ \ if and only if  for any \ $i \in[0, k-h]$ \ the rank \ $h$ \ sub-quotient \ $ F_{i+h}(E)\big/F_{i}(E)$ \ is semi-simple. \\
  So \ $\mathcal{F}_1 : = \mathcal{F}(\lambda_1, \dots, \lambda_k)$ \ and \ $\mathcal{F}_{k}$ \ is the subset of isomorphism classes of semi-simple frescos with fundamental invariants \ $\lambda_{1}, \dots, \lambda_{k}$.
  \item This stratification is invariant by any change of variable.
  \end{enumerate}
  \end{thm}

 \parag{proof}  We shall prove the algebraicity of this stratification and properties i)  by induction on \ $h \geq 1$. For \ $h = 1$ \ the assertions are clear. Assume \ $ h \geq 2$ and that we have already proved algebraicity and assertions i)  for \ $l \leq h-1$. For each \ $i \in [0,k-h]$ \ we have  algebraic maps
 $$ \pi_{i} : \mathcal{F}_{h-1}(\lambda_{1}, \dots, \lambda_{k}) \to \mathcal{F}_{h-1}(\lambda_{i+1}, \dots, \lambda_{i+h}) $$
 given by \ $ [E] \mapsto [F_{i+h}(E)\big/F_{i}(E)]$. Now for each \ $i \in [0,k-h]$ \ we have a polynomial on \ $\mathcal{F}_{h-1}(\lambda_{i+1}, \dots, \lambda_{i+h})$ given by the \ $\beta-$invariant constructed in the proposition \ref{tech.4}. Composed with \ $\pi_{i}$ \ it gives a polynomial \ $ \beta_{i} :  \mathcal{F}_{h-1}(\lambda_{1}, \dots, \lambda_{k}) \to \C $, and, thanks to the induction hypothesis, we have
 $$ \mathcal{F}_{h}(\lambda_{1}, \dots, \lambda_{k}) = \cap_{i \in  [0,k-h+1] } \beta_{i}^{-1}(0) \subset \mathcal{F}_{h-1}(\lambda_{1}, \dots, \lambda_{k}) .$$
 So we obtain the algebraicity of \ $ \mathcal{F}_{h}(\lambda_{1}, \dots, \lambda_{k})$ \ defined by the condition i).\\
 The invariance by change of variable (see [B.11]) is then a direct consequence of the invariance of the principal J.H. sequence with change of variable and the fact that  a change of variable preserves the semi-simplicity. $\hfill \blacksquare$\\
 
 \parag{Remark} The polynomial \ $ \beta : \mathcal{F}_{k-1}(\lambda_{1}, \dots, \lambda_{k}) \to \C $ \ is  quasi-invariant by changes of variable. In the \ $[\lambda]-$primitive case, its weight is equal to \ $w : = p_{1}+ \dots+ p_{k-1}$. This means explicitly that, for any change of variable \ $\theta \in a.\C[[a]]$ \ we have,
  with $ \theta'(0) = \chi(\theta) \not= 0$  \ \ \ $\beta(\theta_{*}(E)) = \chi(\theta)^{w}.\beta(E)$.$\hfill \square$

 \section*{Bibliography}

\begin{itemize}

\item{[Br.70]} Brieskorn, E. {\it Die Monodromie der Isolierten Singularit{\"a}ten von Hyperfl{\"a}chen}, Manuscripta Math. 2 (1970), p. 103-161.


\item{[B.93]} Barlet, D. {\it Th\'eorie des (a,b)-modules I}, in Complex Analysis and Geo-metry, Plenum Press, (1993), p. 1-43.

\item{[B.95]} Barlet, D. {\it Th\'eorie des (a,b)-modules II. Extensions}, in Complex Analysis and Geometry, Pitman Research Notes in Mathematics Series 366 Longman (1997), p. 19-59.


\item{[B.I]} Barlet, D. {\it Sur certaines singularit\'es non isol\'ees d'hypersurfaces I}, Bull. Soc. math. France 134 (2), ( 2006), p.173-200.

\item{[B.II]} Barlet, D. {\it Sur certaines singularit\'es non isolŽes d'hypersurfaces II}, J. Alg. Geom. 17 (2008), p. 199-254.

\item{[B.09]} Barlet, D. {\it P\'eriodes \'evanescentes et (a,b)-modules monog\`enes}, Bollettino U.M.I. (9) II (2009), p. 651-697.

\item{[B.III]} Barlet, D. {\it Sur les fonctions \`a lieu singulier de dimension 1 }, Bull. Soc. math. France 137 (4), (2009),  p. 587-612.

\item{[B.10]} Barlet, D. {\it The theme of a vanishing period}.\\
 English version in  math.arXiv: 1110.1353 (AG-CV).

\item{[B.11]} Barlet, D {\it Construction of quasi-invariant holomorphic parameters for the Gauss-Manin connection of a holomorphic map to a curve}.\\
Preprint  Institut E. Cartan (Nancy) (2012) \ $n^0 7  $, p.1-69 \\
 also in  math.arXiv: 1203.5990 (AG-CV).
.
\item{[B.12]} Barlet,D. {\it A finiteness theorem for \ $S-$relative formal Brieskorn modules}, Preprint  Institut E. Cartan (Nancy) (2012) \ $n^0 12  $, p.1-23 \\
 also in math.arXiv: 1207.4013 (AG-CV).

\item{[B.-S. 04]} Barlet, D. et Saito, M. {\it Brieskorn modules and Gauss-Manin systems for non isolated hypersurface singularities,} J. Lond. Math. Soc. (2) 76 (2007) $n^01$ \ p. 211-224.

\item{[Bj.93]} Bj{\"o}rk, J-E, {\it   Analytic D-modules and applications}, Kluwer Academic  publishers (1993).

\item{[G.65]} Grothendieck, A. {\it On the de Rham cohomology of algebraic varieties} Publ. Math. IHES 29 (1966), p. 93-101.

 \item{[K.76]} Kashiwara, M. {\it b-function and holonomic systems}, Inv. Math. 38 (1976) p. 33-53.

\item{[M.74]} Malgrange, B. {\it Int\'egrale asymptotique et monodromie}, Ann. Sc. Ec. Norm. Sup. 7 (1974), p. 405-430.

\item{[M.75]} Malgrange, B. {\it Le polyn\^ome de Bernstein d'une singularit\'e isol\'ee}, in Lect. Notes in Math. 459, Springer (1975), p. 98-119.

\item{[S.89]} Saito, M. {\it On the structure of Brieskorn lattices}, Ann. Inst. Fourier 39 (1989), p. 27-72.

\end{itemize}

\end{document}